\documentclass[reqno,a4paper,final]{amsart}
\usepackage{amssymb,amstext,amsthm,eucal}
\usepackage{bbold}
\usepackage{array}
\usepackage[latin1]{inputenc}
\newcommand{\ecke}{\;_-\!\rule{0.2mm}{0.2cm}\;\;}

\newcommand{\SO}{\mathrm{SO}}
   
\newcommand{\CSO}{\mathrm{CSO}}

\newcommand{\SU}{\mathrm{SU}}
\newcommand{\U}{\mathrm{U}}

\newcommand{\RR}{\mathbb{R}}
\newcommand{\CC}{\mathbb{C}}

\newcommand{\eL}{\mathcal{L}}
\newcommand{\eB}{\mathcal{B}}

\newcommand{\eN}{\mathcal{N}}
\newcommand{\eR}{\mathcal{R}}
\newcommand{\eS}{\mathcal{S}}
\newcommand{\eT}{\mathcal{T}}

\newcommand{\eQ}{\mathcal{Q}}
\newcommand{\eW}{\mathcal{W}}
\newcommand{\eF}{\mathcal{F}}
\newcommand{\eC}{\mathcal{C}}

\newcommand{\eP}{\mathcal{P}}

\newtheorem{THEO}{\bf Theorem}
\newtheorem{DEF}{\bf Definition}
\newtheorem{PR}{\bf Proposition}

\newtheorem{COR}{\bf Corollary}
\newtheorem{LM}{\bf Lemma}

\newcommand{\MUNCH}[1]{\relax}
\allowdisplaybreaks[1]
\begin{document}
\begin{sloppypar}
\title{About complex structures in conformal tractor calculus}
\author{Felipe Leitner}
\address{Institut f{\"u}r Geometrie und Topologie, Universit{\"a}t Stuttgart, Pfaffenwaldring 57,
D-70569 Stuttgart, Germany}
\email{leitner@mathematik.uni-stuttgart.de}
\date{October 2005}

\begin{abstract} 
The aim of this paper is to describe the geometry of conformal structures in Lorentzian signature, which admit a 
lightlike conformal Killing vector field whose corresponding adjoint tractor acts as complex structure on the standard 
tractor bundle of conformal geometry. Key to the treatment of this problem is CR-geometry and the 
Fefferman construction. In fact, we will consider here partially integrable CR-structures and a slightly generalised 
Fefferman construction for these, which we call the $\ell$-Fefferman construction. We show that a certain class 
of $\ell$-Fefferman metrics on partially integrable CR-spaces provides all solutions to our problem
concerning the complex structures. 
\end{abstract}

\maketitle

\section{Introduction}
\label{ab1} 

A well known construction by Ch. Fefferman (cf. \cite{Fef76}, \cite{BDS77}) assigns to any CR-space an invariantly 
defined 
conformal structure on a circle bundle. This construction provides a method of investigating CR-invariants in the 
realm of conformal geometry. The classical Fefferman construction is restricted to the case of integrable, strictly 
pseudoconvex CR-spaces. However, it is also known that the Fefferman construction applies to more general situations, 
like that of partially integrable CR-geometry, and beyond that to various other parabolic geometries (cf. 
\cite{Cap02}, \cite{Cap05b}). Partially integrable CR-spaces are those which admit a totally real Levi-form and a 
non-trivial Nijenhuis 
torsion tensor.

In this paper we aim to prove a result about conformal geometry in Lorentzian signature, i.e., the metrics in a given 
conformal class have signature $(1,n)$. The problem that we pose admits a reasonably natural formulation only in 
conformal tractor calculus. To solve the problem, we will essentially use the Fefferman construction. However, we need 
to do this for partially integrable CR-spaces. Moreover, we will extent the Fefferman construction in this 
case slightly to something, which we call the $\ell$-Fefferman construction.

To explain our problem, let us consider the adjoint tractor bundle $\mathcal{A}(F)$ on a space $(F,c)$ with conformal 
structure $c$ of Lorentzian signature. Via the so-called splitting operator $\eS$, any conformal vector field $V$ on 
$F$ corresponds to a uniquely and invariantly defined section $\eS(V)$ in $\mathcal{A}(F)$ which solves the 
conformally covariant 
tractor equation \[ \nabla^{nor}\eS(V)=-\Omega^{nor}(V,\cdot)\ , \] where $\nabla^{nor}$ denotes the unique normal 
tractor connection of conformal geometry and $\Omega^{nor}$ is the corresponding curvature $2$-form. On the other 
side, 
any adjoint tractor which solves such an equation stems from a conformal vector field on $F$
(cf. \cite{Cap05c}). The adjoint tractor 
bundle is generated by the adjoint action of the M{\"o}bius group $\SO(2,n+1)$ on $\frak{so}(2,n+1)$ and sections of 
$\mathcal{A}(F)$ can be naturally interpreted as bundle endomorphisms on the standard tractor bundle $\eT(F)$, which 
comes from the standard representation of $\SO(2,n+1)$ on $\RR^{2,n+1}$. We ask the following question in this paper: 
Which Lorentzian conformal structures admit a conformal Killing vector $V$ such that the corresponding adjoint tractor 
$\eR=\eS(V)$ (which necessarily solves the above equation) acts as complex structure on $\eT(F)$?

The answer to this question is known in case that the curvature expression $\Omega^{nor}(V,\cdot)$ on the right 
hand side of the tractor equation is zero. Namely, all spaces providing a solution to this simplified equation are 
exactly the Fefferman spaces of strictly pseudoconvex and integrable CR-spaces
(cf. \cite{Spa85}, \cite{Gra87}). A natural idea to extend the 
description of spaces with solutions is to apply the Fefferman construction for partially integrable CR-spaces.
The latter spaces are known to admit canonical normal Cartan 
connections (cf. \cite{CS00}). The Nijenhuis torsion tensor should then be the (only) contributor to the conformal 
curvature term on the 
right hand side of above tractor equation for $\eS(V)$. We will show that this extension of the classical case is  
one part of the solution to our problem concerning the complex structures on $\eT(F)$. In fact,
it turns out that there is still a 
gap in the geometric description. This is due to the fact that in the curvature term $\Omega^{nor}(V,\cdot)$ 
on the right hand side of the 
tractor equation there can be involved another term which should be seen as a correction to the Weyl connection form 
used in the Fefferman construction. This effect will be incorporated in our $\ell$-Fefferman construction.
    
We describe shortly the course of our investigations. We will start by introducing partially integrable CR-spaces. To 
pursue the Fefferman construction we will use pseudo-Hermitian structures and certain preferred linear connections of 
those. In particular, we will extend the notion of Tanaka-Webster connections to the case of partially integrable 
CR-geometry (cf. \cite{Lee86}, \cite{Miz93}). With the help of the preferred connections we can construct Fefferman 
metrics, whose 
conformal classes 
turn out to be CR-invariants. An important part of our investigation is then the calculation of the relation between 
Webster scalar curvature and Riemannian scalar curvature in the Fefferman construction (cf. Theorem \ref{TH2}). 
Furthermore, we are able to find an explicit expression for the Laplacian applied to the fundamental Killing vector 
field in the Fefferman construction (cf. Proposition \ref{PR1}). It turns out that all the calculations can be 
conducted for the more general $\ell$-Fefferman metrics without further expenses. The results identify the explicit 
form of the adjoint tractor 
that belongs to the fundamental Killing vector in the $\ell$-Fefferman construction (cf. Proposition \ref{PR3}). 
A method of reconstruction shows that the $\ell$-Fefferman construction provides all possible solutions to 
our problem concerning the complex structures on $\eT(F)$ (cf. Proposition \ref{PR5}).
Finally, we summarise our result in Theorem \ref{TH3}.   

\section{CR-manifolds and pseudo-Hermitian geometry}
\label{ab2}

We recall in this section the definition
of CR-structures on smooth manifolds and explain certain notions about their 
integrability conditions. The basic integrability condition that we assume throughout this paper will be the partial 
integrability.  For convenience, we will introduce two equivalent definitions for CR-structures, a complex and real 
version, and use both of them in the sequel. Moreover, we will introduce pseudo-Hermitian forms on CR-manifolds. All 
these concepts are known and for further informations about this subject we refer to e.g. \cite{Lee86}, 
\cite{Bau99} and \cite{Cap02}.

To start with, let $M^n$ be a connected smooth manifold of odd dimension $n=2m+1$. 
First, we introduce the notion of a complex almost
CR-structure on $M$, which is by definition a complex subbundle $T_{10}$ of the complexified tangent 
bundle $TM^{\CC}=TM\otimes\CC$ such that
\[ T_{10}\cap\overline{T_{10}}=\{0\}\quad\mbox{and}\quad dim_{\CC}T_{10}=m\ .\]
We set $T_{01}:=\overline{T_{10}}$ and denote by $\Gamma(T_{10})$ the space of smooth sections
in $T_{10}$ over $M$.
All (complex) almost CR-structures that we will consider shall be 
non-degenerate.
To express this condition, we consider the Levi-form $L$ on $T_{10}$, which is defined by
\[\begin{array}{crcl}L: &T_{10}\times T_{10} &\to& E:=TM^{\CC}/_{T_{10}\oplus T_{01}}\ ,\\[2mm]
&L(U,V)&:=&i\cdot(pr_E[U,\bar{V}])\end{array}\]
where $pr_E$ denotes the projection of (complex) vectors onto the quotient $E$. The complex almost CR-structure 
$T_{10}$ is 
called non-degenerate if its Levi-form $L$ is non-degenerate. Furthermore, 
we want to state certain notions of integrability 
for $T_{10}$. We say that a non-degenerate complex CR-structure $T_{10}$ is
 \[\begin{array}{ccl} \mbox{partially\ integrable} \qquad&\mbox{iff}&\quad [\Gamma(T_{10}),\Gamma(T_{10})]\subset
\Gamma(T_{10}\oplus T_{01})\qquad \mbox{and}\\[2mm]
\mbox{integrable}\quad&\mbox{iff}&\quad [\Gamma(T_{10}),\Gamma(T_{10})]\subset
\Gamma(T_{10})\ .\end{array}\]

The real version of an almost CR-structure consists of a real subbundle $H$ of codimension $1$
in the tangent bundle $TM$ and an almost complex structure $J$ on $H$, i.e., $J^2=-id|_H$. 
The Lie bracket of vector fields on $M$ induces the tensorial map
\[\begin{array}{crcl}L_H: &H\times H &\to& Q:=TM/_{H}\ ,\\[2mm]
&L_H(X,Y)&:=&pr_Q[X,Y]\end{array}\]
where $X,Y\in \Gamma(H)$ are arbitrary vector fields.
The non-degeneracy condition for a real almost CR-structure says
that $H$ is a contact distribution in $TM$, or equivalently, the map 
$L_H$ is non-degenerate.
The partial integrability for $(H,J)$ means that the map $L_H$ is totally real, i.e.,
\[L_H(X,Y)=L_H(JX,JY)\qquad\mbox{for\ all}\ X,Y\in H\ .\]
Finally, the integrability of $(H,J)$ is determined by the additional vanishing of the Nijenhuis tensor:
\[\eN_J(X,Y):=[X,Y]-[JX,JY]+J[JX,Y]+J[X,JY]=0\] for all $X,Y\in\Gamma(H)$. 

The natural correspondence between complex and real version of almost CR-structures is given in the one direction 
starting with $T_{10}$ by 
\[H:=Re(T_{10}\oplus T_{01})\quad\mbox{and}\quad J(U+\bar{U}):=i(U-\bar{U}),\ \  U\in\Gamma(T_{10})\ .\]
In the other direction starting with $(H,J)$, the eigenspaces of the extended complex linear map $J$ on $H^\CC$   
to the eigenvalue $i$ define an almost complex CR-structure $T_{10}$. The introduced notions of integrability 
and non-degeneracy coincide under the natural correspondence.
In this paper, we will be concerned with non-degenerate, partially integrable CR-structures.
We call a pair $(M,T_{10})$ (resp. a triple $(M,H,J)\ $) a partially integrable CR-manifold
if $T_{10}$ (resp. $(H,J)\ $) is non-degenerate and partially integrable.
Often we will not distinguish between the real and the complex version but use both notions
simultaneously. 

A nowhere vanishing $1$-form $\theta$ on a non-degenerate almost CR-manifold $(M,H,J)$ with $\theta|_H\equiv 0$
is called a pseudo-Hermitian form (resp. pseudo-Hermitian structure) and 
the data $(M,H,J,\theta)$ (resp. $(M,T_{10},\theta)$) are called a 
pseudo-Hermitian manifold. 
A pseudo-Hermitian form is necessarily a contact form on $M$.
The pseudo-Hermitian form $\theta$ determines uniquely the Reeb vector field $T\in\frak{X}(M)$ by  
\[\theta(T)\equiv 1\qquad\mbox{and}\qquad d\theta(T,\cdot)\equiv 0\ .\]
The Hermitian form 
\[\begin{array}{crcc}L_\theta: &T_{10}\times T_{10} &\to& \CC\ ,\\[2mm]
&L_\theta(U,V)&:=&-id\theta(U,\bar{V})\end{array}\]
is called the Levi-form of $(M,T_{10},\theta)$. Obviously, it holds $\theta(L(U,V))=L_\theta(U,V)$.
The Levi-form $L_\theta$ can be naturally extended to $TM^\CC$ by
\[\begin{array}{l} L_\theta(U,\bar{V}):=0\ ,\qquad L_\theta(T,\cdot)=0\ ,\\[2mm]
L_\theta(\bar{U},\bar{V}):=\overline{L_\theta(U,V)}=L_\theta(V,U)\ .\end{array}\]
The real part of the extension of $L_\theta$ is a symmetric bilinear form on $TM$, which
is non-degenerate on $H$. We denote this form also by
\[L_\theta: TM\times TM \to \RR\ .\]

In general, the signature of $L_\theta$ on $H$ is $(2p,2q)$ for 
some non-negative integers $p,q$ with $p+q=m$, where  $2p$ is the 
number 
of timelike vectors in an orthonormal frame of $(H,L_\theta)$. Since two 
pseudo-Hermitian forms on $(M,H,J)$ differ only by rescaling with a nowhere vanishing function, the definiteness
of $L_\theta$ is an invariant of the almost CR-structure. Hence, in the definite case, we say that
the almost CR-space $(M,H,J)$ 
is strictly pseudoconvex and 
a pseudo-Hermitian form on such a space is assumed to be positive definite. We will consider in this paper
strictly pseudoconvex spaces. However, this restriction is not essential for our considerations.

\section{The Tanaka-Webster connection} 
\label{ab3}

We assume from now on that $(M,T_{10})$ is a partially integrable, strictly pseudoconvex CR-manifold equipped with 
pseudo-Hermitian 
form $\theta$, i.e., the Levi-form $L_\theta$ is positive definite on $H$. The purpose of this section
is to introduce a certain covariant derivative which naturally belongs to the given pseudo-Hermitian structure.
We call it the Tanaka-Webster connection, since it generalises the classical case for such connections (cf. 
\cite{Lee86}, \cite{Miz93}). 
We also introduce notions of curvature for this connection. 

The following Lemma \ref{LM1} to \ref{LM4} are known facts, certainly for the case of integrable CR-structures, where 
its statements and proofs can be found in (cf. \cite{Bau97}, \cite{Bau99}). We give here corresponding modified 
statements for the 
weaker 
condition of partial integrability.  Thereby, we mainly explain where refinements of the formulae 
concerning  the integrable 
case have to be taken into consideration due to partial integrability. We will not repeat the parts of the proofs 
which do not depend on the Nijenhuis tensor.

\begin{LM} \label{LM1} (cf. \cite{Bau99}) Let $L_\theta:TM^\CC\times TM^\CC\to \CC$ be the Levi-form to $\theta$ on a 
partially
integrable CR-manifold $(M,T_{10})$ (resp. $(M,H,J)$) and let $T$ be the Reeb vector field 
belonging to $\theta$. Then
\[
\begin{array}{l}
[T,Z]\in \Gamma(H^\CC)\qquad \mbox{for\ all}\ Z\in\Gamma(H^\CC)\ ,\\[2mm]
L_\theta([T,U],V)+L_\theta(U,[T,V])=T(L_\theta(U,V))\ ,\\[2mm]
L_\theta([T,\bar{U}],V)=L_\theta([T,\bar{V}],U)\ ,\\[2mm]
L_\theta([T,\U],\bar{V})=L_\theta([T,V],\bar{U})
\end{array}\]
for all $U,V\in\Gamma(T_{10})$, and  
\[
\begin{array}{l}
L_\theta(X,Y)=d\theta(X,JY)\ ,\\[2mm]
L_\theta(JX,JY)=L_\theta(X,Y)\ ,\qquad L_\theta(JX,Y)+L_\theta(X,JY)=0\ ,\\[2mm]
L_\theta([T,X],Y)-L_\theta([T,Y],X)=L_\theta([T,JX],JY)-L_\theta([T,JY],JX)
\end{array}
\]
for all $X,Y\in \Gamma(H)$.
\end{LM}

As next, we state the existence of a particular covariant derivative, which we call the Tanaka-Webster 
connection,
belonging to any pseudo-Hermitian structure $\theta$ on a partially integrable CR-manifold $(M,T_{10})$ (cf. 
\cite{Lee86}, \cite{Bau99}, \cite{Miz93}). 
 
\begin{LM} \label{LM2} (cf. \cite{Bau99}) Let $(M,T_{10},\theta)$ be a pseudo-Hermitian manifold and 
let $T$ be the Reeb vector field to $\theta$. Then there exists a uniquely determined
covariant derivative 
\[\nabla^W:\Gamma(T_{10})\to \Gamma(T^*M^\CC\otimes T_{10})\]
such that
\begin{enumerate}
\item
\[\nabla^W_TU=pr_{10}[T,U]\ , \qquad \nabla^W_{\bar{V}}U=pr_{10}[\bar{V},U]\ ,\]
\item
\[X(L_\theta(U,V))=L_\theta(\nabla^W_XU,V)+L_\theta(U,\nabla^W_{\bar{X}}V)\]
\end{enumerate}
for all $U,V\in\Gamma(T_{10})$ and $X\in TM^\mathbb{C}$, where $pr_{10}$ denotes the projection onto $T_{10}$.
Moreover, $\nabla^W$ 
satisfies
\[\nabla^W_UV-\nabla^W_VU=pr_{10}[U,V]\ .\]
\end{LM}

We note that in case of an integrable CR-structure $T_{10}$ the more special relation
\[\nabla^W_UV-\nabla^W_VU=[U,V]\]
holds.
Nevertheless, the proof of Lemma \ref{LM2} mainly uses 
Lemma \ref{LM1} and 
remains the same as in 
\cite{Bau97}. We extend now the Tanaka-Webster connection $\nabla^W$ to the complex 
tangent bundle $TM^\CC$ by
\[\nabla^WT:=0\qquad\mbox{and}\qquad\nabla^W_X\bar{U}
:=\overline{\nabla^W_{\bar{X}}U}\qquad \mbox{for\ all}\ X\in TM^{\mathbb{C}},\ U\in\Gamma(T_{10})\ .\]
The torsion $Tor^W$ of this connection is defined in the usual manner as
\[Tor^W(X,Y):=\nabla^W_XY-\nabla^W_YX-[X,Y]\]
for $X,Y\in\Gamma(TM^{\CC})$. 

\begin{LM} \label{LM3} (cf. \cite{Bau99}) The torsion $Tor^W$ of the Tanaka-Webster connection 
\[\nabla^W\ :\ \Gamma(TM^\CC)\to \Gamma(T^*M^\CC\otimes TM^\CC)\]
satisfies
\[\begin{array}{l} Tor^W(U,\bar{V})=iL_\theta(U,V)\cdot T\ ,\\[2mm]
Tor^W(U,V)=-pr_{01}[U,V]\ , \qquad
Tor^W(\bar{U},\bar{V})=-pr_{10}[\bar{U},\bar{V}]\ ,\\[2mm]
Tor^W(T,U)=-pr_{01}[T,U]\ ,\qquad
Tor^W(T,\bar{U})=-pr_{10}[T,\bar{U}]
\end{array}\]
for all $U,V\in \Gamma(T_{10})$, where $pr_{01}$ denotes the projection onto $T_{01}$.
\end{LM}

In the integrable case the formulae for the torsion simplify to 
\[Tor^W(U,V)=Tor^W(\bar{U},\bar{V})=0\] 
when $U,V\in T_{10}$. However, by the refinement in Lemma \ref{LM2}, it holds in the partially integrable case  
only $\nabla^W_UV-\nabla^W_VU-[U,V]=-pr_{01}[U,V]$. This implies the torsion formulae in Lemma \ref{LM3}.
Finally, we can restrict the Tanaka-Webster connection to 
its real part and obtain a linear connection $\nabla^W$ on the (real) tangent bundle $TM$. 

\begin{LM} \label{LM4} (cf. \cite{Bau99}) Let $\theta$ be a pseudo-Hermitian structure on a partially integrable 
CR-manifold
$(M,H,J)$. The Webster connection  
\[\nabla^W:\frak{X}(M)\to\Gamma(T^*M\otimes TM)\]
is uniquely determined by the following properties 
\begin{enumerate}
\item
\[
X(L_\theta(Y,Z))=L_\theta(\nabla^W_XY,Z)+L_\theta(Y,\nabla^W_XZ)\]
for all $X,Y,Z\in\frak{X}(M)$, 
i.e., $\nabla^W$ is metric with respect to $L_\theta$ on $H$,
\item
\[
Tor^W(X,Y)=L_\theta(JX,Y)\cdot T-\frac{1}{4}\eN_J(X,Y)\quad\qquad \mbox{and}\]
\item
\[
Tor^W(T,X)=-\frac{1}{2}\big(\ [T,X]+J[T,JX]\ \big)
\]
\end{enumerate}
for all $X,Y\in\Gamma(H)$. In addition, the connection $\nabla^W$ satisfies
\[\nabla^WT=0\qquad\mbox{and}\qquad\nabla^WJ=0\ .\]
\end{LM} 

{\bf Proof.} Again, the same proof as in \cite{Bau97} (p.16) applies. We have only to take into consideration
the refined formula $Tor^W(U,V)=-pr_{01}[U,V]$ and the
relation
\[ -\frac{1}{4}\eN_J(X,Y)=Tor^W(U,V)+Tor^W(\bar{U},\bar{V})\]
for all $U,V\in T_{10}$, 
whereby $X=U+\bar{U}$ and $Y=V+\bar{V}$. Together with $Tor^W(U,\bar{V})=iL_\theta(U,V)\cdot T$ this shows
\begin{eqnarray*}
Tor^W(X,Y)&=&Tor^W(\bar{U},V)+Tor^W(U,\bar{V})+Tor^W(U,V)+Tor^W(\bar{U},\bar{V})\\
&=&L_\theta(JX,Y)\cdot T-\frac{1}{4}\eN_J(X,Y)
\end{eqnarray*}
for all $X,Y\in H$. Since 
\begin{eqnarray*}
\nabla^W_YJX&=&\nabla^W_Y(i(U-\bar{U}))=i(\nabla_Y^WU-\overline{\nabla_{Y}^WU})\\
&=&J(\nabla^W_YU\ +\ \overline{\nabla^W_YU})=J\nabla^W_YX
\end{eqnarray*}
for all $X\in\Gamma(H)$ and $Y\in TM$, it follows that $J$ is parallel with respect to $\nabla^W$.
\hfill$\Box$\\

The important point in Lemma \ref{LM4} is that the Nijenhuis tensor is part of  
the torsion of the Tanaka-Webster connection due to partial
integrability. 
The way in which the Nijenhuis tensor occurs in the torsion $Tor^W$ is essentially implicated
by condition (1) in Lemma \ref{LM2}. (The metric condition (2) of Lemma \ref{LM2}
is a rather inevitable choice.)
There are other suitable connections, which occur naturally
in the framework of CR-geometry resp. pseudo-Hermitian geometry, most notably the so-called Weyl connections
in the sense of \cite{CS03}. This type of connections differs from the Tanaka-Webster connections
just by its torsion normalisation. We will later meet the Weyl connections on line bundles. 

Finally, we define here expressions of curvature for $\nabla^W$.
Thereby, we introduce the following conventions and notations concerning calculations
with respect to a local basis or frame. The indices with letters $i,j$ and $k$ run from $1$ to $2m$,
whereas the indices with Greek letters $\alpha,\beta$ and $\gamma$ run from $1$ to $m$.
With $\{e_i:\ i=1,\ldots, 2m\}$ we denote
an orthonormal basis (resp. local frame) of $(H,L_\theta)$ such that the additional condition
\[J(e_{2\alpha-1})=e_{2\alpha},\qquad\ \ J(e_{2\alpha})=-e_{2\alpha-1}\qquad\ \mbox{for\ all\ } \alpha=1,\ldots, m\]
is satisfied.
Then the complex vectors \[Z_\alpha:=\frac{1}{\sqrt{2}}(e_{2\alpha-1}-iJe_{2\alpha-1}),\quad \alpha=1,\ldots, m\ ,\]
form an orthonormal basis
of $(T_{10},L_\theta)$ and  the vectors $Z_{\bar{\alpha}}:=\overline{Z_\alpha}$ with  $\alpha=1,\ldots, m$ represent
an orthonormal basis of $T_{01}$. 
The curvature operator of $\nabla^W$ is defined by
\[ R^{\nabla^W}(X,Y)Z=\nabla^W_X\nabla^W_YZ-\nabla^W_Y\nabla^W_XZ-\nabla^W_{[X,Y]}Z\ ,
\]
where $X,Y,Z\in \Gamma(TM^\CC)$ are complex vectors. The curvature operator is tensorial
in $X,Y$ and $Z$. Moreover, we have the curvature tensor given by
\[ R^W(X,Y,Z,V)=L_\theta(R^{\nabla^W}(X,Y)Z,\bar{V})
\]
for $X,Y,Z,V\in TM^\CC$. A straightforward calculation proves that
\begin{eqnarray*}
R^W(X,Y,Z,V)&=& -R^W(Y,X,Z,V)\ =-R^W(X,Y,V,Z),\\[2mm]
R^W(A,\overline{B},C,\overline{D})&=&\ R^W(C,\overline{B},A,\overline{D})\
-\ L_\theta(Tor^W(\overline{B},Tor^W(C,A)),D),\\[2mm]
R^{\nabla^W}(A,B)\overline{C}&=&\ (\nabla_{\overline{C}}Tor^W)(A,B)
\end{eqnarray*}
for all vectors $X,Y,Z,V$ in $TM^\CC$ and $A,B,C,D$ in $T_{10}$.
In case that the Nijenhuis tensor $\eN$ vanishes the two latter identities simplify to
\begin{eqnarray*}
R^W(A,\overline{B},C,\overline{D})&=&\ R^W(C,\overline{B},A,\overline{D})\qquad\mbox{and}\\[2mm]
R^{\nabla^W}(A,B)\overline{C}&=&\ 0\ .
\end{eqnarray*}

The Webster-Ricci and scalar curvatures are defined as contractions of $R^W$:
\[\begin{array}{l}Ric^W:=\sum_{\alpha=1}^mR^W(\cdot,\cdot,Z_\alpha,Z_{\bar{\alpha}})\ ,\\[3mm]
scal^W:=\sum_{\alpha=1}^mRic^W(Z_\alpha,Z_{\bar{\alpha}})\ .\end{array}\]
These definitions are independent of the choice of orthonormal basis. With respect to the $e_i$'s we have
\[\begin{array}{l}Ric^W(X,Y)=i\cdot\sum_{\alpha=1}^m R^W(X,Y,e_{2\alpha-1},Je_{2\alpha-1})\ ,\\[3mm]
scal^W=i\cdot\sum_{\alpha=1}^m Ric^W(e_{2\alpha-1},Je_{2\alpha-1})\ .\end{array}\]
Obviously, the function $scal^W$ on $(M,T_{10},\theta)$ is real.
We set
\[\omega_\alpha^\beta:=L_\theta(\nabla^WZ_\alpha,Z_\beta)\ .\]
With these components of the Tanaka-Webster connection the Ricci and scalar curvature
are expressed by
\begin{eqnarray*}
Ric^W(X,Y)&=&\sum_{\alpha=1}^m L_\theta([\nabla_X^W,\nabla_Y^W]Z_\alpha-\nabla^W_{[X,Y]}Z_\alpha,Z_\alpha)\\[3mm]
&=&\left(\sum_{\alpha=1}^m
d\omega_{\alpha}^{\alpha}-\sum_{\alpha,\beta=1}^m\omega_{\alpha}^{\beta}\wedge\omega_{\beta}^{\alpha}\right)(X,Y)\\[3mm]
&=&\sum_{\alpha=1}^m d\omega_{\alpha}^{\alpha}(X,Y)\ ,\\[4.5mm]
scal^W&=&\sum_{\alpha,\beta=1}^m d\omega_{\alpha}^{\alpha}(Z_\beta,Z_{\bar{\beta}})\ .
\end{eqnarray*}
Note that the Webster-Ricci curvature $Ric^W$ is not symmetric, in general. However, it is
symmetric in case of integrable CR-structures.

\section{Rescaling of a pseudo-Hermitian structure} 
\label{ab4}

We discuss here the transformation rules for the Tanaka-Webster connection and its scalar curvature
under rescaling of a given pseudo-Hermitian structure. We remark right at the beginning  that calculations 
and resulting expressions in this 
section do not differ formally from those in the classical integrable case, since the Nijenhuis torsion does not play 
any 
role in the transformation (cf. \cite{Lee86}). The transformation rule will help us to prove the 
invariance of the Fefferman construction, which happens in the next section.

Let $(M,T_{10},\theta)$ be a partially integrable, strictly pseudoconvex CR-space with pseudo-Hermitian
structure $\theta$ and corresponding Tanaka-Webster connection $\nabla^W$. 
As always, we use local frames $\{e_i:\ i=1,\ldots, 2m\}$ with the property 
\[J(e_{2\alpha-1})=e_{2\alpha},\qquad\ \ J(e_{2\alpha})=-e_{2\alpha-1}\qquad\ \mbox{for\ all\ } \alpha=1,\ldots, m,\]
and we set $Z_\alpha=\frac{1}{\sqrt{2}}(e_{2\alpha-1}-iJe_{2\alpha-1})$. Moreover, we set
\[\theta^{\alpha}=L_\theta(\cdot,Z_{\alpha})\ ,\qquad \theta^{\bar{\alpha}}=L_\theta(\cdot,Z_{\bar{\alpha}}) \]
and
\[\delta_\alpha^\beta=
L_\theta(Z_\alpha,Z_\beta)\ ,\qquad\omega_{\alpha}^{\beta}:=L_{\theta}(\nabla^WZ_\alpha,Z_\beta)\ .\]
Then, it is $\nabla^WZ_\alpha=\sum_{\beta=1}^m\omega_\alpha^\beta\otimes Z_\beta$.
Furthermore, for a real smooth function $f\in C^\infty(M)$ we define
\[\begin{array}{l}
f_\alpha:=Z_\alpha(f)\ ,\qquad f_{\bar{\alpha}}:=Z_{\bar{\alpha}}(f)\ ,\qquad f_o=T(f),\\[3mm]
f_{\alpha\bar{\beta}}:=(\nabla^W_{Z_{\bar{\beta}}}df)(Z_\alpha)\qquad\mbox{and}\qquad
f_{\bar{\alpha}\beta}:=(\nabla^W_{Z_{\beta}}df)(Z_{\bar{\alpha}})\ .
\end{array}\] 
Eventually, we set
\[
\delta f:=\sum_{\alpha=1}^m f_{\bar{\alpha}}Z_\alpha\qquad\mbox{and}\qquad
\Delta_b f:=-\sum_{\alpha=1}^m \big(f_{\alpha\bar{\alpha}}+f_{\bar{\alpha}\alpha}\big)\ .\]
The latter definitions are independent of the choice of orthonormal frame. It holds 
$\delta f(f)=\sum_{\alpha=1}^m f_\alpha\cdot f_{\bar{\alpha}}$ 
and the differential operator $\Delta_b$ is called sublaplacian.

Now let $f\in C^\infty(M)$ be an arbitrary function and let $\tilde{\theta}=e^{2f}\theta$.
The Hermitian form $L_{\tilde{\theta}}$ is again positive definite and $\tilde{\theta}$ is 
another pseudo-Hermitian structure on $(M,T_{10})$.
In fact, any two pseudo-Hermitian forms on $(M,T_{10})$ differ only by multiplication with a smooth 
positive function. 
We want to examine the transformation rules for the Tanaka-Webster connections and the corresponding
scalar curvatures under such a rescaling.
First of all, we notice that
\[\tilde{\theta}^\alpha=e^f(\theta^\alpha+2if_{\bar{\alpha}}\theta)\qquad \mbox{for\ all\ }\ \alpha=1,\ldots,m, \]
where $\tilde{\theta}^\alpha$ is dual to 
$\tilde{Z}_{\alpha}=e^{-f}\cdot Z_{\alpha}$ with respect to $L_{\tilde{\theta}}$,
and
\[ \tilde{T}=e^{-2f}\cdot (T-2i\delta f+2i\overline{\delta f})\ .\]

\begin{LM} \label{LM5} (cf. \cite{Lee86}) Let $\tilde{\theta}=e^{2f}\theta$ be a rescaled pseudo-Hermitian structure 
on 
a partially integrable CR-space $(M,T_{10})$. Then it holds 
\begin{enumerate}
\item
\begin{eqnarray*} \tilde{\nabla}^W_XU\quad  
=& &\nabla^W_XU+2df(U)\cdot pr_{10}X+2df(pr_{10}X)\cdot U-2L_{\theta}(U,\bar{X})\cdot\delta f\\[3mm]  
&+&i\theta(X)\cdot \left( 4df(U)\cdot\delta f+4U\cdot \delta f(f)+2\cdot\nabla^W_U\delta f\right)\ \end{eqnarray*}
for all $U\in\Gamma(T_{10})$ and $X\in TM^\CC$.
\item  For the connection components it holds 
\begin{eqnarray*} \tilde{\omega}_\alpha^\beta 
&=& \omega_\alpha^\beta\ +\ 2(f_\alpha\theta^\beta-f_{\bar{\beta}}\theta^{\bar{\alpha}})
+\delta_\alpha^\beta\cdot\left(\sum_{\gamma=1}^m 
f_\gamma\theta^\gamma-f_{\bar{\gamma}}\theta^{\bar{\gamma}}\right)\\[1.5mm] 
&&\quad\ \ \ +\ i\ \theta(\cdot) \left(
f_{\bar{\beta}\alpha}+f_{\alpha\bar{\beta}}+4f_\alpha f_{\bar{\beta}}+4\delta_\alpha^\beta\cdot
\sum_{\gamma=1}^m f_\gamma 
f_{\bar{\gamma}}\right)\ ,\end{eqnarray*}
where $\ \alpha,\beta= 1,\ldots, m$.
\item
The Webster scalar curvature rescales by
\[\widetilde{scal}^W=e^{-2f}\cdot \left(\ scal^W\ +\ 2(m+1)\Delta_bf\ -\ 4m(m+1)\delta f(f)\ \right)\quad .
\]
\end{enumerate}
\end{LM}

{\bf Proof.} (1) We take the expression for $\tilde{\nabla}^W$ in Lemma \ref{LM5} as definition 
for a connection and verify that it satisfies the
determining properties for the Tanaka-Webster connection of Lemma \ref{LM2}. 
First, we have
\[ \tilde{\nabla}^W_{\bar{V}}U=\nabla^W_{\bar{V}}U-2L_\theta(U,V) \delta f\ ,\]
and on the other side, it is 
\[\tilde{pr}_{10}[\bar{V},U]=pr_{10}[\bar{V},U]-2L_\theta(U,V)\delta f\]
for all $U,V\in\Gamma(T_{10})$. 
This shows that $\tilde{\nabla}^W_{\bar{V}}U=\tilde{pr}_{10}[\bar{V},U]$.
Next we see that
\begin{eqnarray*}
&&\quad L_{\tilde{\theta}}(\tilde{\nabla}^W_XU,V)+L_{\tilde{\theta}}(U,\tilde{\nabla}^W_{\bar{X}}V)\\[3mm]&=&
\quad  e^{2f}\cdot L_\theta(\nabla^W_XU+2df(U)pr_{10}X+2df(pr_{10}X)U-2L_\theta(U,\bar{X})\delta f,V)\\
&&+\ e^{2f}\cdot L_\theta(U,\nabla^W_{\bar{X}}V+2df(V)pr_{10}\bar{X}+2df(pr_{10}\bar{X})V-2L_\theta(V,X)\delta f)\\[3mm]
&=&\quad e^{2f}\cdot X(L_\theta(U,V))+2e^{2f}df(X)\cdot L_\theta(U,V)\\
&& +\ e^{2f}\cdot \left( 2df(U)\cdot L_\theta(X,V)-2df(\bar{V})\cdot L_\theta(U,\bar{X})\right)\\
&&+\ e^{2f}\cdot \left( 2df(\bar{V})\cdot L_\theta(U,\bar{X})-2df(U)\cdot \overline{L_\theta(V,X)}\right)\\[3mm]
&=& X(L_{\tilde{\theta}}(U,V))\end{eqnarray*}
for all $U,V\in\Gamma(T_{10})$. It remains to show that $\tilde{\nabla}^W_{\tilde{T}}U=
\tilde{pr}_{10}[\tilde{T},U]$. For this purpose we calculate 
\begin{eqnarray*} e^{2f}\cdot \tilde{pr}_{10}[\tilde{T},U]&=&\quad pr_{10}[T,U]+2ie^{2f}pr_{10}
[e^{-2f}(-\delta f+\overline{\delta f}),U]\\[3mm]
&=&\quad pr_{10}[T,U]-4idf(U)\delta f-2i\cdot \nabla^W_{\delta f}U+2i\cdot \nabla^W_{\overline{\delta f}}U\\
&&+\ 2i\cdot \sum_{\alpha=1}^m U(f_{\bar{\alpha}})\cdot 
Z_\alpha+2i\cdot\sum_{\alpha=1}^m f_{\bar{\alpha}}\cdot\nabla^W_UZ_\alpha\ ,
\end{eqnarray*}
and on the other side, 
\begin{eqnarray*} e^{2f}\cdot \tilde{\nabla}^W_{\tilde{T}}U&=&\quad\nabla^W_TU-2i\nabla^W_{\delta f}U
+2i\nabla^W_{\overline{\delta f}}U-8i\cdot df(U)\delta f-4i\sum_{\alpha=1}^m f_{\alpha}f_{\bar{\alpha}}\cdot U\\
&&+\ i(4df(U)\delta f+4\delta f(f)\cdot U+2\cdot \nabla^W_U\delta f)\ , 
\end{eqnarray*}
which obviously equals the previous expression.

(2) To calculate the connection components $\tilde{\omega}_\alpha^\beta$, we use the
following identity
\begin{eqnarray*}&&if_o\delta_\alpha^\beta + 2\cdot L_\theta(\nabla_{Z_\alpha}^W\delta f,Z_\beta)\\[2.5mm]
&=& d\theta(Z_\alpha,Z_{\bar{\beta}})T(f)+2\cdot Z_\alpha(f_{\bar{\beta}})-2\cdot \nabla^W_{Z_{\bar{\alpha}}}Z_\beta(f)\\
&=& 2\cdot Z_\alpha(f_{\bar{\beta}})+[Z_{\bar{\beta}},Z_\alpha](f)-pr_{10}[Z_{\bar{\beta}},Z_\alpha](f)
-pr_{01}[Z_\alpha,Z_{\bar{\beta}}](f)\\
&=&Z_\alpha(f_{\bar{\beta}})+Z_{\bar{\beta}}(f_\alpha)-\nabla^W_{Z_\alpha}Z_{\bar{\beta}}(f)
-\nabla^W_{Z_{\bar{\beta}}}Z_\alpha(f)\\[2.5mm]
&=& f_{\bar{\beta}\alpha}+f_{\alpha\bar{\beta}}\ .\end{eqnarray*}
It is
\begin{eqnarray*} \tilde{\omega}_\alpha^\beta&=&\quad 
L_{\tilde{\theta}}(\tilde{\nabla}^W\tilde{Z}_\alpha,\tilde{Z}_\beta)
= L_\theta(\tilde{\nabla}^WZ_\alpha,Z_\beta)-\delta_\alpha^\beta\cdot df\\[3mm]
&=&\quad \omega_\alpha^\beta-\delta_\alpha^\beta\cdot df+2\cdot 
f_\alpha\theta^\beta-2f_{\bar{\beta}}\theta^{\bar{\alpha}}+2\delta_\alpha^\beta\cdot 
\sum_{\gamma=1}^m f_\gamma\theta^\gamma\\[-1mm]
&&+\ i\cdot\left( 4f_\alpha f_{\bar{\beta}}
+2\cdot L_\theta(\nabla^W_{Z_\alpha}\delta f,Z_\beta)+
4\delta_\alpha^\beta\cdot \sum_{\gamma=1}^m f_\gamma f_{\bar{\gamma}}
\right)\theta\\[3mm]
&=&\quad 
\omega_\alpha^\beta+2(f_\alpha\theta^\beta-f_{\bar{\beta}}\theta^{\bar{\alpha}})+\delta_\alpha^\beta\cdot
\sum_{\gamma=1}^m
(f_\gamma\theta^\gamma-f_{\bar{\gamma}}\theta^{\bar{\gamma}})\\[-1mm]
&&+\ i\cdot\left( 4f_\alpha f_{\bar{\beta}}+if_o\delta_\alpha^\beta
+2\cdot L_\theta(\nabla^W_{Z_\alpha}\delta f,Z_\beta)+
4\delta_\alpha^\beta\cdot \sum_{\gamma=1}^m f_\gamma f_{\bar{\gamma}}
\right)\theta\\[3mm]
&=&\quad \omega_\alpha^\beta + 
2(f_\alpha\theta^\beta-f_{\bar{\beta}}\theta^{\bar{\alpha}})+\delta_\alpha^\beta\cdot\sum_{\gamma=1}^m
(f_\gamma\theta^\gamma-f_{\bar{\gamma}}\theta^{\bar{\gamma}})\\
&&+\ i\cdot\left( 4f_\alpha f_{\bar{\beta}}+
f_{\bar{\beta}\alpha}+f_{\alpha\bar{\beta}}+4\delta_\alpha^\beta\cdot\sum_{\gamma=1}^m f_\gamma 
f_{\bar{\gamma}}\right)\theta\quad . \end{eqnarray*}

(3) We use now the latter formula for the connection components to calculate the Webster scalar curvature.
It is
\begin{eqnarray*} \sum_{\alpha=1}^m \tilde{\omega}_\alpha^\alpha &=\quad&\quad\sum_{\alpha=1}^m \omega_\alpha^\alpha\ 
\ +\ 
\ 
(m+2)\cdot\sum_{\alpha=1}^m(f_\alpha\theta^\alpha-f_{\bar{\alpha}}\theta^{\bar{\alpha}})\\
&\quad+& i\cdot\sum_{\alpha=1}^m \left(\ f_{\bar{\alpha}\alpha}+f_{\alpha\bar{\alpha}}\  +\  4(m+1)f_\alpha 
f_{\bar{\alpha}}\right)\theta\quad .
\end{eqnarray*}
The trace of the exterior differential of this expression 
is equal to the Webster scalar curvature. We use
\begin{eqnarray*}d(\sum_{\alpha=1}^m f_\alpha\theta^\alpha)(Z_\gamma,Z_{\bar{\beta}})&=& 
\sum_{\alpha=1}^m (df_\alpha\wedge \theta^\alpha +f_\alpha d\theta^\alpha)(Z_\gamma,Z_{\bar{\beta}})\\
&=&
-Z_{\bar{\beta}}(f_\gamma)
+\sum_{\alpha=1}^m f_\alpha\theta^\alpha([Z_{\bar{\beta}},Z_\gamma])\\
&=& -f_{\gamma\bar{\beta}}
\end{eqnarray*}
and obtain
\begin{eqnarray*}\sum_{\alpha,\beta=1}^m d\tilde{\omega}_\alpha^\alpha(Z_\beta,Z_{\bar{\beta}})&=\quad&\ \
\sum_{\alpha,\beta=1}^m d\omega_\alpha^\alpha(Z_\beta,Z_{\bar{\beta}})\ \ 
-\ \ (m+2)\cdot\sum_{\beta=1}^m (f_{\beta\bar{\beta}}
+f_{\bar{\beta}\beta})\\
&\quad+& i\sum_{\alpha,\beta=1}^m\left(\ f_{\bar{\alpha}\alpha}+f_{\alpha\bar{\alpha}}\  +\  4(m+1)\cdot f_\alpha
f_{\bar{\alpha}}\ \right)\cdot d\theta(Z_\beta,Z_{\bar{\beta}})\ \ .\end{eqnarray*}
As result we have
\[\widetilde{scal}^W=e^{-2f}\cdot\left(\ scal^W\ +\ 2(m+1)\Delta_b f\ -\ 4m(m+1)\delta f(f) \ \right)\ \ \ . \]
{ }\hfill $\Box$

\section{The Fefferman metric} 
\label{ab5}

We construct now the Fefferman metric to a pseudo-Hermitian structure 
$\theta$ on the total space of the canonical $S^1$-principal bundle of
a partially integrable, strictly pseudoconvex CR-manifold
$(M,T_{10})$. The importance of the Fefferman metric
relies on the independence of its conformal class from the chosen pseudo-Hermitian
structure, i.e., the Fefferman conformal class is an invariant of
the underlying CR-structure. The construction that we describe 
coincides with the classical Fefferman construction for integrable CR-spaces.
In fact, it looks formally the same as in the classical case (cf. \cite{Lee86}, \cite{Bau99}).
However, at the end of the section we aim to introduce a slightly more generalised class 
of metrics, which we call the $\ell$-Fefferman metrics. 
The reason for this extension shall find its justification in
the later sections.  
 
Let $(M^n,T_{10})$ be a partially integrable, strictly pseudoconvex CR-manifold of dimension $n=2m+1$ 
and let $\theta$
be a pseudo-Hermitian structure on this CR-space.  We denote by
\[
\Lambda^{m+1,0}M\ := \ \{\ \rho \in \Lambda^{m+1}M\otimes\CC :\ X\ecke \rho=0 \ 
\mbox{\ for \ all\ } X\in T_{01}=\overline{T_{10}}
\ \}\] 
the complex line bundle over $M^n$, which consists of all those  
complex $(m+1)$-forms that vanish when an element of $T_{01}$ is inserted.
The bundle $\Lambda^{m+1,0}M$ is called the canonical line bundle of the CR-space $(M, T_{10})$. 
The positive real numbers $\RR^+$ act by multiplication on $K^*:=\Lambda^{m+1,0}\setminus \{0\}$, which denotes the
canonical line bundle without zero section.
We set $F_c:=K^*/_{\RR^+}$ and the triple 
\[(\ F_c,\pi,M\ )\] denotes the canonical $S^1$-principal bundle of $(M,T_{10})$
whose fibre action is induced by complex multiplication with the elements of the unit circle $S^1$ in $\CC$. 

Let $\{Z_\alpha:\ \alpha=1,\ldots, m\}$ be a local orthonormal frame of $(T_{10},L_\theta)$
and let $\theta^\alpha,$\ $\alpha=1,\ldots,m$, denote the corresponding dual $1$-forms. The $(m+1,0)$-form 
\[
\tau:=\theta\wedge\theta^1\wedge\ldots\wedge\theta^m
\]
is a local section of $\Lambda^{m+1,0}M$. We denote by $[\tau]$ the corresponding
local section in $F_c=K^*/_{\RR^+}$. With help of the projection $\pi$ every $1$-form $\rho$ on $M$ can be lifted
to $F_c$. The result is a $1$-form $\pi^*\rho$ on $F_c$. For convenience, we shall usually denote 
the lifted $1$-form by $\rho$ again. With help of the section $[\tau]$ we can also pull back
$1$-forms $\sigma$ to $M$, which is denoted by $[\tau]^*\sigma$ or 
just $\sigma$ again. 

The Tanaka-Webster connection $\nabla^W$ naturally extends to a covariant derivative on sections of the 
$(m+1,0)$-form bundle  $\Lambda^{m+1,0}(M)$. In fact, the covariant derivative $\nabla^W$ acting
on $\Lambda^{m+1,0}M$ is induced by a unique connection $1$-form on the $S^1$-principal fibre bundle
$F_c$, which we denote by
\[A^W: TF_c\ \to\ i\RR\ \ .\]  
It is $[\tau]^*A^W=-\sum_{\alpha=1}^m\omega_{\alpha}^\alpha$ with respect to the local frame forms 
$\theta^\alpha, \alpha=1,\ldots,m$.
Further, we set
\[A_\theta:= A^W-\frac{i}{2(m+1)} scal^W\theta\ .\]
The latter is a connection $1$-form on $F_c$ as well. It is the so-called Weyl 
connection on the canonical $S^1$-bundle of 
$(M,T_{10})$, which 
belongs to the given pseudo-Hermitian structure $\theta$ (cf. \cite{CS03}).
The curvature of $A^W$ is the $2$-form $\Omega^W=dA^W$. 
It holds 
\[\Omega^W=-\sum_{\alpha=1}^m d\omega_\alpha^\alpha=-Ric^W\ .\] 
We denote the curvature of $A_\theta$ by $\Omega_\theta=dA_\theta$. It is
\[\Omega_\theta=-Ric^W-\frac{i}{2(m+1)}scal^W d\theta-\frac{i}{2(m+1)}d(scal^W)\cdot\theta\ .\]

We define now the Fefferman metric to $\theta$ on $F_c$ by
\[ f_\theta\ :=\ \pi^*L_\theta-i\frac{4}{m+2}\pi^*\theta\circ A_\theta\ ,\]
or in shorter notation, we often use the expression $f_\theta=L_\theta-i\frac{4}{m+2}\theta\circ A_\theta$.
This is, in fact, a symmetric $2$-tensor on the real tangent bundle of $F_c$. 
In case of an underlying  strictly pseudoconvex space its 
signature is Lorentzian (i.e., $(1,2m+1)$).
The Fefferman conformal class $[f_\theta]$ 
consists of all smooth metrics $\tilde{f_\theta}$ on $F_c$ which arise by conformal rescaling of $f_\theta$,
i.e., it is $\tilde{f_\theta}=e^{2\phi}f_\theta$ for some smooth function $\phi$.
As we already mentioned the Fefferman conformal class shall be independent of the particular choice
of a pseudo-Hermitian structure on $M$. This makes the Fefferman conformal class $[f_\theta]$ an 
invariant object in a natural manner attached to the CR-structure $T_{10}$ on $M$. We want to prove this
invariance property for $f_\theta$.
For this purpose we need to find the transformation rule for the Tanaka-Webster connection form $A^W$ resp. for 
the Weyl connection $A_\theta$ on $F_c$ under rescaling of $\theta$. 
Thereby, we use the results from the last section.  
We will see that the transformation rule for the Weyl connection is particular
easy.        

To start with, let $\tilde{\theta}=e^{2f}\cdot\theta$ be a rescaled pseudo-Hermitian structure on $(M,T_{10})$
and let $\tilde{Z}_\alpha=e^{-f}\cdot Z_\alpha$,  where $\alpha=1,\ldots,m$, be the rescaled basis vectors. It is 
\[\tilde{\theta}^\alpha=e^f(\theta^\alpha+2if_{\bar{\alpha}}\theta)
\qquad\mbox{and}\qquad\tilde{\tau}=e^{(m+2)f}\tau\ ,\] 
i.e., $[\tilde{\tau}]$ and $[\tau]$
are identical as local sections of $F_c$. Moreover, it is
\[\begin{array}{l} [\tau]^*A^W=-\sum_{\alpha=1}^m\omega_\alpha^\alpha\quad\mbox{and}\\[4mm]
 {}[\tau]^*\tilde{A}^W=-\sum_{\alpha=1}^m\tilde{\omega}_\alpha^\alpha\ .\end{array}\]
By using Lemma \ref{LM5}, we calculate
\begin{eqnarray*}
[\tau]^*(\tilde{A}^W-A^W)&=&-\ (m+2)\cdot\sum_{\alpha=1}^m (f_\alpha\theta^\alpha
-f_{\bar{\alpha}}\theta^{\bar{\alpha}})\\ 
&&+\ i(\ \Delta_b f\ -\ 4(m+1)\cdot\sum_{\alpha=1}^mf_\alpha f_{\bar{\alpha}}\ )\cdot \theta,\end{eqnarray*}  
and further,
\begin{eqnarray*}[\tau]^*(A_{\tilde{\theta}}-A_\theta)&=& -\ (m+2)\cdot\sum_{\alpha=1}^m (\ f_\alpha\theta^\alpha-
f_{\bar{\alpha}}\theta^{\bar{\alpha}}\ )\\
&&-\ i(\ -\Delta_bf\ +\ 4m(m+1)\cdot\sum_{\alpha=1}^m f_\alpha f_{\bar{\alpha}}\ )\cdot \theta\\
&&-\ \frac{i}{2(m+1)}\left(scal^W\!+ 2(m+1)\Delta_b f-4(m+1)\cdot\sum_{\alpha=1}^m f_\alpha f_{\bar{\alpha}}\right)
\!\cdot\theta\\
&&+\ \frac{i}{2(m+1)}\cdot scal^W\theta\\[4mm]
&=&-\ (m+2)\cdot\sum_{\alpha=1}^m(\ f_\alpha\theta^\alpha-f_{\bar{\alpha}}\theta^{\bar{\alpha}}\ )\\
&&-\ \left(2i(m+2)\cdot\sum_{\alpha=1}^m f_\alpha f_{\bar{\alpha}}\right)\cdot\theta\qquad .
\end{eqnarray*}
We conclude that
\[A_{\tilde{\theta}}=A_\theta\ -\ (m+2)\sum_{\alpha=1}^m(\ 
f_\alpha\theta^\alpha-f_{\bar{\alpha}}\theta^{\bar{\alpha}}\ )\
-\ 2i(m+2)\sum_{\alpha=1}^mf_\alpha f_{\bar{\alpha}}\theta\ .\]
Now we can consider the transformation rule for the Fefferman metric $f_\theta$ under rescaling
of $\theta$. It is
\begin{eqnarray*} f_{\tilde{\theta}}&=&\quad 2\cdot \sum_{\alpha=1}^m 
\tilde{\theta}^\alpha\circ\tilde{\theta}^{\bar{\alpha}}\
\ \ -\ \ \ i\frac{4}{m+2}\tilde{\theta}\circ A_{\tilde{\theta}}\\
&=&\quad e^{2f}\cdot\big(\ \sum_{\alpha=1}^m2\cdot\big(\ \theta^\alpha\circ\theta^{\bar{\alpha}}
+2if_{\bar{\alpha}}\theta\circ
\theta^{\bar{\alpha}}-2if_\alpha\theta^{\alpha}\circ\theta+4f_{\bar{\alpha}}f_\alpha\cdot \theta\circ\theta\ \big)\\
&&\qquad -\ i\frac{4}{m+2}\theta\circ A_\theta +\sum_{\alpha=1}^m\big(\ 
4if_\alpha\theta\circ\theta^\alpha-4if_{\bar{\alpha}}\theta
\circ\theta^{\bar{\alpha}}-8f_\alpha f_{\bar{\alpha}}\cdot\theta\circ\theta\ \big)\ \big)\\
&=&\quad e^{2f}\cdot f_\theta\ \ .\
\end{eqnarray*}

\begin{THEO} \label{TH1} (cf. \cite{Lee86}) Let $(M^n,T_{10})$ be a partially integrable CR-space with a 
pseudo-Hermitian 
structure $\theta$ and Fefferman metric $f_\theta$ on $F_c$. 
Let $\tilde{\theta}=e^{2f}\theta$ be a rescaled pseudo-Hermitian structure. Then the 
corresponding
Fefferman metric rescales by $f_{\tilde{\theta}}=e^{2f}\cdot f_\theta$.
\end{THEO}

As we have defined Fefferman metrics in this section they 
depend on the connection $A_\theta$ associated to $\theta$. 
A more general class of metrics on $F_c$ (which are again Lorentzian for strictly 
pseudoconvex CR-structures) is given by choosing any connection 
form $A$. The difference between $A$ and $A_\theta$ is the
lift of a $1$-form $\ell$ on $M$ with purely imaginary values. We denote
\[
A_{\theta,\ell}:=A_\theta+\ell\qquad\mbox{and}\qquad 
f_{\theta,\ell}=L_\theta-i\frac{4}{m+2}\theta\circ A_{\theta,\ell}\ .
\]
\begin{DEF} Let $(M^n,T_{10})$ be a partially integrable CR-space,
$\theta$ a pseudo-Hermitian structure and $\ell\in\Omega^1(M;i\RR)$ an arbitrary 
$1$-form on $M$. Then we call the metric
\[f_{\theta,\ell}=L_\theta-i\frac{4}{m+2}\theta\circ A_{\theta,\ell}
\]
the (generalised) $\ell$-Fefferman metric with respect to $\theta$
on $(M^n,T_{10})$.
\end{DEF}
In case that $\ell=0$ the metric $f_{\theta,\ell}=f_{\theta}$ is just the
usual Fefferman metric to $\theta$. If $\ell$ is a closed form on $M$
then the metrics $f_\theta$ and $f_{\theta,\ell}$ are locally isometric. 
The local
isometry is just given by a gauge transformation on $F_c$ which transforms the 
connection form $A_{\theta,\ell}$ into $A_\theta$. In particular,
the fibres are preserved. On the other side, a local isometry
between $f_\theta$ and $f_{\theta,\ell}$, which preserves the fibres, can only exist if there is a gauge 
transformation,
i.e., if the difference $\ell=A_{\theta,\ell}-A_\theta$ is closed. Moreover, metrics
$f_{\theta,\ell}$ and $f_{\tilde{\theta},\tilde{\ell}}$ are locally isometric
only if the rescaling function $f$ is constant zero, i.e., when $\tilde{\theta}=\theta$ holds. Altogether,
we can conclude that generalised Fefferman metrics $f_{\theta,\ell}$ and $f_{\tilde{\theta},\tilde{\ell}}$ 
are locally isometric with preserved fibre if and only if $\tilde{\theta}=\theta$ and $\tilde{\ell}-\ell$ 
is a closed form on $M$. Finally, we notice that if $\tilde{\theta}=e^{2f}\theta$ and 
$\ell\in\Omega^1(M;i\RR)$ is arbitrary then 
\[f_{\tilde{\theta},\ell}=e^{2f}\cdot f_{\theta,\ell}\ .\]
This shows that the conformal class $[f_{\theta,\ell}]$ of the
$\ell$-Fefferman metric is an invariant of the pair $(T_{10},\ell)$,
which consists of a partially integrable CR-structure and a $1$-form $\ell\in\Omega^1(M;i\RR)$.
We denote the conformal class of the $\ell$-Fefferman metric by $c_\ell:=[f_{\theta,\ell}]$,
where $\theta$ is some arbitrary pseudo-Hermitian form.

\section{The torsion tensor}    
\label{ab6}

We examine here properties of the torsion tensor $Tor^W$ with respect to the Tanaka-Webster connection on 
pseudo-Hermitian spaces. The torsion consists essentially of the Nijenhuis tensor $\eN_J$ and $Tor^W(T,\cdot)$, where 
the latter part is the deviation from transversal symmetry along $T$. The Nijenhuis tensor is a CR-invariant. The 
discussion in this section will be for use 
in the next section when the torsion enters our calculations for curvature expressions 
that arise in connection with the Fefferman construction.

Let $(M^n,H,J)$ be a strictly pseudoconvex, partially integrable CR-space with dimension $n=2m+1$
and let $\theta$
denote a pseudo-Hermitian structure on $M$. By $\{e_i:\ i=1,\ldots,2m\}$ we denote
an orthonormal basis of $L_\theta$ on $H$ such that
\[ J(e_{2\alpha-1})=e_{2\alpha}\qquad\mbox\qquad J(e_{2\alpha})=-e_{2\alpha-1}\qquad \mbox{for\ all\ \ }
\alpha=1,\ldots,m\ . \]
Now we introduce the following conventions.
Let $A$ be a $(2,r)$-tensor (with $r=0$ or $1$). We denote its trace (or contraction) by
\[tr_\theta A:=\sum_{i=1}^{2m}\ A(e_i,e_i)\ .\]
More generally, we use the notation
$tr_\theta^{k,l}A$ for the trace of a $(s,r)$-tensor $A$, where the contraction
takes place in the $k$th and $l$th entry of $A$. In relation with the contraction we also use
the convention 
$(\nabla A)(X,Y)=(\nabla_XA)(Y)$, where
$\nabla$ denotes some
covariant derivative on
the tangent space.
If $A$ is a skew-symmetric tensor on $M$ then we set
\[L_A(\cdot,\cdot):=iA(\cdot,J\cdot)\ .\] It is
\[tr_\theta L_A=i\cdot tr_\theta A(\cdot,J\cdot)=2\cdot\sum_{\alpha=1}^m A(Z_\alpha,Z_{\bar{\alpha}})\ .\]
If $A$ is a symmetric $(2,0)$-tensor then \[tr_\theta A=2\cdot \sum_{\alpha=1}^m A(Z_\alpha,Z_{\bar{\alpha}})\ .\]

The Nijenhuis tensor $\eN$ on $M$ is defined as
\[\eN_J(X,Y):=[X,Y]-[JX,JY]+J[JX,Y]+J[X,JY]\ ,\]
whereby $X,Y\in \Gamma(H)$. Since $Tor^W(X,Y)=L_\theta(JX,Y)\cdot T-\frac{1}{4}\eN(X,Y)$,
the Nijenhuis tensor can be considered as the essential part of the torsion
restricted to the contact distribution $H$.
It holds
\begin{eqnarray*}
J\eN(X,Y)\ \ =\ \ -\eN(JX,Y)&=&-\eN(X,JY)\qquad\mbox{and}\\[2mm]
tr_\theta L_\theta(\eN(X,\cdot),\cdot)&=&\quad 0
\end{eqnarray*}
for all $X,Y,Z\in H$.
We form with $\eN$
the $\eB_\theta$-tensor by
\[\eB_\theta(X,Y,Z):=\frac{1}{8}\big(\ L_\theta(\eN(X,Y),Z) + L_\theta(\eN(Z,Y),X)
+L_\theta(\eN(Z,X),Y)\ \big)\ .
\]
Moreover,  let \[\eB(X,Y):=\sum_{i=1}^{2m}\eB_\theta(X,Y,e_i)e_i\]
denote the corresponding $(2,1)$-tensor.
The tensor $\eB$ does not depend on the chosen $\theta$ and
the orthonormal frame $\{e_i\}$.
It holds
\[\eB(X,Y)-\eB(Y,X)=\frac{1}{2}\eN(X,Y)\]
and $\eB$ vanishes identically if and only if $\eN$ vanishes identically.
In other words, the tensor $\eB$ contains the same information as $\eN$.
Moreover, it holds
\[\begin{array}{rcl}
\eB_\theta(X,Y,Z)&=&-\eB_\theta(X,Z,Y)\ ,\\[3mm]
\eB_\theta(X,Y,Z)&=&
-\eB_\theta(JX,JZ,Y)=-\eB_\theta(JX,Y,JZ)=
-\eB_\theta(X,JY,JZ)\ ,\\[3mm]
J\eB(X,Y)&=&-\eB(JX,Y)=-\eB(X,JY)\ ,\\[3mm]
tr_\theta\eB&=&\quad\sum_{i=1}^{2m}\eB(e_i,e_i)=0\ ,\\[3mm]
tr_\theta^{1,3}B_\theta(X)&=&\quad\sum_{i=1}^{2m}\eB_\theta(e_i,X,e_i)=0\ ,\\[3mm]
tr_\theta^{2,3}B_\theta(X)&=&\quad\sum_{i=1}^{2m}\eB_\theta(X,e_i,e_i)=0\ .
\end{array}\]
Straightforward calculations using essentially the condition of
partial integrability show the following identities:
\begin{eqnarray*}
tr_\theta L_\theta(\eN(\eN(X,\cdot),Y),\cdot)&=&
tr_\theta L_\theta(\eN(\eN(X,\cdot),\cdot),Y)
-tr_\theta L_\theta(\eN(X,\cdot),\eN(Y,\cdot))\\[2mm]
tr_\theta B_\theta(\eN(X,\cdot),\cdot,Y)&=&\frac{1}{4} tr_\theta L_\theta(\eN(X,\cdot),\eN(Y,\cdot))\\[2mm]
tr_\theta L_\theta(\eB(X,\cdot),\eB(Y,\cdot))&=&\frac{1}{8} tr_\theta
L_\theta(\eN(\eN(X,\cdot),\cdot),Y)\\[2mm]
tr_\theta L_\theta(\eB(X,\cdot),\eB(\cdot,Y))&=&\frac{1}{16} tr_\theta
L_\theta(\eN(\eN(X,\cdot),\cdot),Y)\\[2mm]
\sum_{i,j=1}^{2m}L_\theta(\eN(\eN(e_i,e_j),e_j),e_i)&=&
\frac{1}{2}\cdot\sum_{i,j=1}^{2m}L_\theta(\eN(e_i,e_j),\eN(e_i,e_j))\ .
\end{eqnarray*}
The third identity above shows  that the tensor $tr_\theta L_\theta(\eN(\eN(X,\cdot),\cdot),Y)$ is symmetric in 
$X$
and
$Y$.

The second part of the torsion is
\[Tor^W(T,X)=-\frac{1}{2}(\ [T,X]+J[T,JX]\ )\ ,\]
where $X\in H$. We define the tensor
\[\eT_\theta(X,Y):=-2\cdot L_\theta(Tor^W(T,X),Y)\qquad\mbox{for}\ \ X,Y\in H\ . \]
Then it holds
\begin{eqnarray*}
Tor^W(T,JX)&=&-J(Tor^W(T,X)),\\[3mm]
\eT_\theta(X,Y)&=&\eT_\theta(Y,X)\ =\ L_\theta([T,X],Y)+L_\theta([T,Y],X),\\[3mm]
\eT_\theta(X,JY)&=&\eT_\theta(JX,Y)\qquad\quad\qquad\mbox{and}\\[3mm]
tr_\theta\eT_\theta&=&
tr_\theta\eT_\theta(\cdot,J\cdot)\ =\ 0\ .
\end{eqnarray*}

\section{The scalar curvature of a Fefferman metric} 
\label{ab7}

We calculate in this section parts of the  Ricci-curvature tensor and with the help of
these the scalar
curvature for Fefferman 
metrics. 
In fact, we will do this calculation for the more general class of $\ell$-Fefferman metrics.
The formulae that we obtain generalise the results of \cite{Lee86}. Due to partial
integrability the Nijenhuis tensor $\eN$ will enter the curvature expressions.
This makes our calculations more laborious.

Let $(M^n,H,J)$ be a strictly pseudoconvex, partially integrable CR-space with dimension $n=2m+1$ and
let $\theta$ denote a pseudo-Hermitian structure on $M$. 
The $\ell$-Fefferman metric on the canonical $S^1$-bundle $F_c$ to $\theta$ on $M$ is defined as 
\[f_{\theta,\ell}=L_\theta-i\frac{4}{m+2}\theta\circ A_{\theta,\ell}\ ,\]
where $\ell$ is (the lift of) an arbitrary $1$-form on $M$ with purely imaginary values and
\[A_{\theta,\ell}:=A_\theta+\ell\ ,\] i.e., $A_{\theta,\ell}$ 
takes the form of a generic
connection $1$-form on $F_c$. In this section, we will sometimes denote the $\ell$-Fefferman
metric abbreviated by $f$.

The $S^1$-action on the fibres of $F_c$ induces
a (vertical) fundamental vector field for each element in the
Lie algebra $i\RR$ of $S^1$. We denote by $S$ the fundamental field which is determined
by
\[ A_{\theta,\ell}(S)=i\frac{m+2}{2}\ \ \ . \]
The field $S$ is lightlike on $(F_c,f_{\theta,\ell})$.
As usual, we denote by  
$\{e_i:\ i=1,\ldots,2m\}$
an orthonormal local frame of
$(H,L_\theta)$ satisfying \[J(e_{2\alpha-1})=e_{2\alpha}\qquad \mbox{and}\qquad
J(e_{2\alpha})=-e_{2\alpha-1}\qquad \mbox{for\ all\ }\alpha=1,\ldots,m\ .\]
Moreover, let $T$ be the Reeb vector field to $\theta$ on $M$ and let
$X^*$ denote the horizontal lift to $F_c$ of any vector $X$ on $M$ with respect to the connection $A_{\theta,\ell}$.
It holds $f_{\theta,\ell}(S,T^*)=1$ and
\[\{e_1^*,\ldots,e_{n}^*,T^*,S\}\]
is a local frame on $(F_c,f_{\theta,\ell})$.
Throughout this section we use
(local) vector fields $X,Y,Z$ and $V$ on $M$ which have constant coefficients
with
respect to the chosen local frame $\{e_i:\ i=1,\ldots,2m\}$. In particular,  
this implies that scalar products of $X^*,Y^*,Z^*$ and $V^*$
with each other and respect to $f_{\theta,\ell}$ are constant.
To start with the calculations, we note that
\[\begin{array}{lcl}
[X^*,S]&=&0\\[3mm]
{}[X^*,Y^*]_{Vert}&=&i\frac{2}{m+2}\Omega_{\theta,\ell}(X,Y)\cdot S\\[3mm]
{}[X^*,Y^*]_{Horiz}&=&[X,Y]^*\\[3mm]
{}[T^*,X^*]&=&[T,X]^*+i\frac{2}{m+2}\Omega_{\theta,\ell}(T,X)\cdot S\\[3mm]
{}[X^*,Y^*]&=&pr_H[X,Y]^*-d\theta(X,Y)\cdot T^*+i\frac{2}{m+2}\Omega_{\theta,\ell}(X,Y)\cdot S\ ,
\end{array}\]
where $\Omega_{\theta,\ell}=dA_{\theta,\ell}$ is the curvature of the connection form $A_{\theta,\ell}$ on $F_c$.
\begin{LM} \label{LM6} (cf. \cite{Bau99}) For the $\ell$-Fefferman metric $f$ on  $F_c$, it holds
\[\begin{array}{ccl}
f(\nabla^f_{X^*}Y^*,Z^*)&=&\quad L_\theta(\nabla^W_XY,Z)+\eB_\theta(X,Y,Z)\\[3mm]
f(\nabla^{f}_{S}Y^*,Z^*)&=&\quad \frac{1}{2}L_\theta(JY,Z)\\[3mm]
f(\nabla^{f}_{X^*}Y^*,S)&=&-\ \frac{1}{2}L_\theta(JX,Y)\\[3mm]
f(\nabla^{f}_{T^*}Y^*,Z^*)&=&\quad\frac{1}{2}\big(\ L_\theta([T,Y],Z)-L_\theta([T,Z],Y)
-i\frac{2}{m+2}\Omega_{\theta,\ell}(Y,Z)\ \big)\\[3mm]
f(\nabla^{f}_{X^*}Y^*,T^*)&=&\quad\frac{1}{2}\big(\ L_\theta([T,X],Y)+L_\theta([T,Y],X)
+i\frac{2}{m+2}\Omega_{\theta,\ell}(X,Y)\ \big)\\[3mm]
f(\nabla^{f}_{T^*}T^*,Z^*)&=&-\ i\frac{2}{m+2}\Omega_{\theta,\ell}(T,Z)\\[3mm]
f(\nabla^{f}S,S)&=&\quad f(\nabla^{f}S,T^*)=f(\nabla^{f}T^*,T^*)=0\\[3mm]
f(\nabla^{f}_SS,Z^*)&=&\quad f(\nabla^{f}_ST^*,Z^*)=f(\nabla^{f}_{T^*}S,Z^*)=0
\end{array}\]
for all $X,Y,Z\in\Gamma(H)$, which have pairwise constant scalar products with respect to $L_\theta$.
\end{LM}

{\bf Proof.} We apply the Koszul formula for the Levi-Civita connection $\nabla^{f}$, namely
it holds
\[f(\nabla^{f}_DB,C)=\frac{1}{2}\big(\ f([D,B],C)+f([C,B],D)
+f([C,D],B)\ \big)\]
for all vector fields $B,C,D$ on $F_c$, which have constant length and pairwise constant
scalar products. In fact, the formulae of Lemma \ref{LM6} result immediately from this Koszul
formula, the above expressions for commutators of vector fields on $F_c$
and replacing the
scalar products with
respect to $f$
by those with respect to $L_\theta$ after projecting the vectors to $M$.
For example, for the first formula we find that
\begin{eqnarray*} 
2f(\nabla^{f}_{X^*}Y^*,Z^*)&=&L_\theta([X,Y],Z)+L_\theta([Z,Y],X)+L_\theta([Z,X],Y)\\[1mm]
&=&2L_\theta(\nabla^W_XY,Z)
+\frac{1}{4}L_\theta(\eN_J(X,Y),Z)+\frac{1}{4}L_\theta(\eN_J(Z,Y),X)\\[1mm]
&&\qquad\qquad\qquad\ \  \!\! +\frac{1}{4}L_\theta(\eN_J(Z,X),Y)\\[1mm]
&=&2L_\theta(\nabla^W_XY,Z)+2\eB_\theta(X,Y,Z)
\end{eqnarray*}
for all sections $X,Y,Z$ in $H$ with pairwise constant scalar products. The other formulae
follow in a similar way by applying the Koszul formula and the above expressions for commutators.
\hfill$\Box$\\

\noindent
Note that in the expressions for the connection components in Lemma \ref{LM6} the Nijenhuis
torsion occurs only in the $H$-part of the horizontal distribution. The $1$-form
$\ell$ influences the curvature expression that appear.

\begin{LM} \label{LM7}(cf. \cite{Lee86}) It is
\begin{eqnarray*}
Ric^{f_{\theta,\ell}}\; (S\; ,T^*)&=\quad\ &\!\!\!\frac{1}{2(m+1)}scal^W\ - 
\ \frac{i}{2(m+2)}tr_\theta d\ell(\cdot,J\cdot)\qquad\quad\mbox{and}\\[6mm]
Ric^{f_{\theta,\ell}}(X^*,V^*)&=\quad\ &\!\!\!\frac{scal^W}{(m+1)(m+2)}\cdot L_\theta(X,V)\\[3.5mm]
&\quad\ +&\!\!\!i\frac{m}{2(m+2)}(\ Ric^W(X,JV)+Ric^W(V,JX)\ )\\[3.5mm]
&\quad\ -&\!\!\!\frac{m}{4}(\ \eT_\theta(X,JV)+\eT_\theta(V,JX)\ )\\[3.5mm]
&\quad\ +&\!\!\!tr_\theta^{1,4}(\nabla^W\eB_\theta)(X,V)+tr_\theta^{1,4}(\nabla^W\eB_\theta)(V,X)\\[3.5mm]
&\quad\ -&\!\!\!\frac{1}{8}tr_\theta L_\theta(\eN(X,\cdot),\eN(V,\cdot))\
+\ \frac{1}{4}tr_\theta L_\theta(\eN(\eN(X,\cdot),\cdot),V)\\[3.5mm]
&\quad\ +&\!\!\!\frac{i}{m+2}\left(\ d\ell(X,JV)+d\ell(V,JX)\ \right)
\end{eqnarray*}
for all vectors $X^*,V^*$ in the horizontal lift of $H$ to $TF_c$.
\end{LM}

{\bf Proof.}
We will use the connection components of Lemma \ref{LM6} in order to obtain 
second covariant derivatives of vector fields on $F_c$ and certain components
of the Riemannian curvature tensor $R^{f_{\theta,\ell}}$. 
For convenience, we set
\[ G(X,V):=\eT_\theta(X,V)+i\frac{2}{m+2}\Omega_{\theta,\ell}(X,V)\]
for all $X,V$ in $H$. 
First, we have
\begin{eqnarray*}
f(\nabla^{f}_{X^*}\nabla^{f}_{Y^*}Z^*,V^*)&=\quad&
\!\!\!X^*(f(\nabla_{Y^*}^{f}Z^*,V^*))
 -f(\nabla^{f}_{Y^*}Z^*,\nabla^{f}_{X^*}V^*)\\[4mm]
&=\quad&\!\!\! L_\theta(\nabla^W_X\nabla^W_YZ,V)\\[1mm]
&\quad+&\!\!\! \frac{1}{4}L_\theta(JY,Z)\cdot G(X,V) + \frac{1}{4}L_\theta(JX,V)\cdot G(Y,Z)\\[1mm]
&\quad-&\!\!\! L_\theta(\nabla^W_YZ,\eB(X,V))+L_\theta(\nabla^W_X(\eB(Y,Z)),V)\\[1mm]
&\quad-&\!\!\! L_\theta(\eB(Y,Z),\eB(X,V))\ ,
\end{eqnarray*}
\begin{eqnarray*}
f(\nabla^{f}_{[X^*,Y^*]}Z^*,V^*)&=\quad&\!\!\!
L_\theta(\nabla^W_{pr_H[X,Y]}Z,V)+\frac{i}{m+2}\Omega_{\theta,\ell}(X,Y)\cdot L_\theta(JZ,V)\\[1mm]
&\quad-&\!\!\!\frac{1}{2}L_\theta(JX,Y)\cdot\\
&&\!\!\!\big(\ L_\theta([T,Z],V)-L_\theta([T,V],Z)-i\frac{2}{m+2}\Omega_{\theta,\ell}(Z,V)\ \big)\\[1mm]
&\quad+&\!\!\!\eB_\theta(pr_H[X,Y],Z,V)\ \ \ ,
\end{eqnarray*}
which results to the curvature component
\begin{eqnarray*}
R^{f}(X^*,Y^*,Z^*,V^*)&=\quad&\!\!\!R^W(X,Y,Z,V)\\[1.5mm]
&\quad-&\!\!\!\frac{i}{m+2}L_\theta(JZ,V)\cdot\Omega_{\theta,\ell}(X,Y)-\frac{1}{2}L_\theta(JX,Y)\cdot G(Z,V)\\[1.5mm]
&\quad-&\!\!\!L_\theta(JX,Y)\cdot L_\theta(Tor^W(T,Z),V)\\[1.5mm]
&\quad+&\!\!\!\frac{1}{4}L_\theta(JY,Z)\cdot G(X,V)+\frac{1}{4}L_\theta(JX,V)\cdot G(Y,Z)\\[1.5mm]
&\quad-&\!\!\!\frac{1}{4}L_\theta(JX,Z)\cdot G(Y,V)-\frac{1}{4}L_\theta(JY,V)\cdot G(X,Z)\\[1.5mm]
&\quad-&\!\!\!(\nabla^W_Y\eB_\theta)(X,Z,V)+(\nabla^W_X\eB_\theta)(Y,Z,V)\\[1.5mm]
&\quad-&\!\!\!\frac{1}{4}\eB_\theta(\eN(X,Y),Z,V)\\[1.5mm]
&\quad-&\!\!\!L_\theta(\eB(Y,Z),\eB(X,V))+L_\theta(\eB(X,Z),\eB(Y,V))
\end{eqnarray*}
for $X,Y,Z$ and $V$ in $H$.
Moreover, it is
\begin{eqnarray*}
f(\nabla^{f}_{T^*}\nabla^{f}_{X^*}S,V^*)\!&=&\!
-f(\nabla^{f}_{X^*}S,\nabla^{f}_{T^*}V^*)\\[1.5mm]
&=&\!-\frac{1}{4}(L_\theta([T,V],JX)-L_\theta([T,JX],V)-\frac{2i}{m+2}\Omega_{\theta,\ell}(V,JX)),\\[4mm]
f(\nabla^f_{[X^*,T^*]}S,V^*)&=&\ \ \frac{1}{2}L_\theta([T,X],JV)
\end{eqnarray*}
and we obtain
\begin{eqnarray*}
R^{f}(X^*,T^*,S,V^*)+R^{f}(X^*,S,T^*,V^*)\!\!\!\!&=\!\!\!\!&-
\frac{i}{2(m+2)}\big(\Omega_{\theta,\ell}(V,JX)+\Omega_{\theta,\ell}(X,JV)\big)\\[2mm]
&&-\frac{1}{4}\left(\ \eT_\theta(X,JV)+\eT_\theta(V,JX)\ \right)\ .
\end{eqnarray*}
With these curvature components we calculate that
\begin{eqnarray*}
\sum_{i=1}^{2m}R^{f_{\theta,\ell}}(X^*,e_i^*,e_i^*,V^*)&=\quad&\!\!\!\sum_{i=1}^{2m}R^W(X,e_i,e_i,V)\\[0.5mm]
&\quad+&\!\!\! \frac{3i}{2(m+2)}(\ \Omega_{\theta,\ell}(X,JV)-\Omega_{\theta,\ell}(JX,V)\ )\\[0.5mm]
&\quad-&\!\!\!\sum_{i=1}^{2m}\big(\ (\nabla^W_{e_i}\eB_\theta)(X,e_i,V)
+\frac{1}{4}\eB_\theta(\eN(X,e_i),e_i,V)\ \big)\\[0.5mm]
&\quad+&\!\!\!\sum_{i=1}^{2m}L_\theta(\eB(X,e_i),\eB(e_i,V))\ ,\end{eqnarray*}
\begin{eqnarray*}
\sum_{i=1}^{2m}R^{f_{\theta,\ell}}(X^*,e_i^*,Je_i^*,JV^*)&=\quad&\!\!\!\sum_{i=1}^{2m}R^W(X,e_i,e_i,V)\\[0.5mm]
&\quad-&\!\!\!\frac{1}{2(m+1)}scal^WL_\theta(X,V)+i\frac{m+1}{m+2}\Omega_{\theta,\ell}(X,JV)\\[1.5mm]
&\quad-&\!\!\!\frac{m-1}{2}L_\theta(Tor^W(T,X),JV)\\[1.5mm]
&\quad-&\!\!\!\frac{m-1}{2}L_\theta(Tor^W(T,V),JX)\\[0.5mm]
&\quad+&\!\!\!\sum_{i=1}^{2m}\big(\ (\nabla^W_{e_i}\eB_\theta)(X,e_i,V)
+\frac{1}{4}\eB_\theta(\eN(X,e_i),e_i,V)\ \big)\\[0.5mm]
&\quad-&\!\!\!\sum_{i=1}^{2m}L_\theta(\eB(X,e_i),\eB(e_i,V))
\end{eqnarray*}
and
\begin{eqnarray*}
\sum_{i=1}^{2m}R^{f_{\theta,\ell}}(X^*,JV^*,e_i^*,Je_i^*)&=\quad&\sum_{i=1}^{2m}R^W(X,JV,e_i,Je_i)\\[0.5mm]
&\quad+&\!\!\!\frac{\ scal^W}{m+1}L_\theta(X,V)\ -\ i\frac{2m+2}{m+2}\Omega_{\theta,\ell}(X,JV)\\[0.5mm]
&\quad-&\!\!\!2\cdot\sum_{i=1}^{2m}L_\theta(\eB(X,e_i),\eB(V,e_i))\ \ \ .
\end{eqnarray*}
By using
\begin{eqnarray*}
\sum_{i=1}^{2m}R^{f_{\theta,\ell}}(X^*,e_i^*,Je_i^*,JV^*)\!\!&=&\sum_{\alpha=1}^{m}R^{f_{\theta,\ell}}
(X^*,e_{2\alpha-1}^*,Je_{2\alpha-1}^*,JV^*)\\
&&\!\!\!\!\!-\sum_{\alpha=1}^{m}R^{f_{\theta,\ell}}(X^*,Je_{2\alpha-1}^*,e_{2\alpha-1}^*,JV^*)\\
&=&\!\!\!\!\!-\sum_{\alpha=1}^{m}R^{f_{\theta,\ell}}(X^*,JV^*,e_{2\alpha-1}^*,Je_{2\alpha-1}^*)
\end{eqnarray*}
we obtain
\begin{eqnarray*}
\sum_{i=1}^{2m}R^{f_{\theta,\ell}}(X^*,e_i^*,e_i^*,V^*)\!\!\!&=\quad&\!\!\!iRic^W(X,JV)\\
&\quad+&\!\!\!\frac{3i}{2(m+2)}(\ \Omega_{\theta,\ell}(X,JV)-\Omega_{\theta,\ell}(JX,V)\ )\\[1.5mm]
&\quad+&\!\!\!\frac{m-1}{2}(L_\theta(Tor^W\!(T,X),JV)\!+\!L_\theta(Tor^W\!(T,V),JX))\\
&\quad-&\!\!\!\sum_{i=1}^{2m}\big(\ 2(\nabla^W_{e_i}\eB_\theta)(X,e_i,V)+\frac{1}{2}\eB_\theta(\eN(X,e_i),e_i,V)\ 
\big)\\
&\quad+&\!\!\!\sum_{i=1}^{2m}\!\big(2L_\theta(\eB(X,e_i),\!\eB(e_i,V))\!+\!L_\theta(\eB(X,e_i),\!\eB(V,e_i))\big).
\end{eqnarray*}
Adding the expression for
$R^{f_\theta}(X^*,T^*,S,V^*)+R^{f_\theta}(X^*,S,T^*,V^*)$ and using the identities for 
torsion terms of section \ref{ab6}
we obtain
\begin{eqnarray*}
Ric^{f_{\theta,\ell}}(X^*,V^*)\!\!\!&=\quad&\!\!\!iRic^W(X,JV)\ +\ 
2\cdot tr_\theta^{1,4}(\nabla^W\eB_\theta)(X,V)\\
&\quad+&\!\!\!\frac{i}{m+2}(\ \Omega_{\theta,\ell}(X,JV)+\Omega_{\theta,\ell}(V,JX)\ )\\[1.5mm]
&\quad-&\!\!\!\frac{m}{4}(\ \eT_\theta(X,JV)+\eT_\theta(V,JX)\ )\\[1.5mm]
&\quad-&\!\!\!\frac{1}{8}tr_\theta L_\theta(\eN(X,\cdot),\eN(V,\cdot))\
+\ \frac{1}{4}tr_\theta L_\theta(\eN(\eN(X,\cdot),\cdot),V)\ .
\end{eqnarray*}
After symmetrisation in $X$ and $V$ of (the first line of) the right hand side of the latter equation,
we obtain the component
$Ric^{f_{\theta,\ell}}(X^*,V^*)$ as stated. Furthermore, it is
\begin{eqnarray*}
Ric^{f_{\theta,\ell}}(S,T^*)\!\!\!&=\quad&\!\!\!\!\!\!\sum_{i=1}^{2m}\big(\ 
f_{\theta,\ell}(\nabla^{f_{\theta,\ell}}_{[T^*,e_i^*]}S,e_i^*)
-f_{\theta,\ell}(\nabla^{f_{\theta,\ell}}_{T^*}\nabla^{f_{\theta,\ell}}_{e_i^*}S,e_i^*)\ \big)\\[1.5mm]
&=\quad&\!\!\!\!\!\!\frac{i}{2(m+2)}\sum_{i=1}^{2m}\big(\  \Omega_{\theta,\ell}(Je_i,e_i)
-L_\theta([T,e_i],Je_i)+L_\theta([T,Je_i],e_i)\ \big)\\[1.5mm]
&=\quad&\!\!\!\!\!\!\frac{1}{2(m+1)}scal^W\ - \ \frac{i}{2(m+2)}tr_\theta d\ell(\cdot,J\cdot)\qquad .
\end{eqnarray*}
{$ \ $} \hfill$\Box$
\begin{THEO} \label{TH2} (cf. \cite{Lee86}) The scalar curvature of the $\ell$-Fefferman metric $f_{\theta,\ell}$  
on $F_c$ is given by 
\[ scal^{f_{\theta,\ell}}=\frac{2m+1}{m+1}\cdot scal^W\ +\ \frac{1}{m+2}\cdot tr_\theta L_{d\ell}\qquad .\]
\end{THEO}
{\bf Proof.}
It is 
\begin{eqnarray*}
scal^{f_{\theta,\ell}}&=\quad&\!\!\!2\cdot Ric^{f_{\theta,\ell}}(T^*,S)+\sum_{i=1}^{2m}
Ric^{f_{\theta,\ell}}(e_i^*,e_i^*)\\[2.5mm]
&=\quad&\!\!\!\left(\ \frac{1}{m+1}+\frac{2m}{(m+1)(m+2)}+\frac{2m}{m+2}\ \right)\cdot scal^W\\[1.5mm]
&\quad -&\!\!\!\frac{1}{8}\cdot\sum_{i,j=1}^{2m}L_\theta(\eN(\eN(e_i,e_j),e_j),e_i)\\[1.5mm]
&\quad +&\!\!\!\frac{1}{4}\cdot\sum_{i,j=1}^{2m}L_\theta(\eN(e_i,e_j),\eN(e_i,e_j))\\[1.5mm]
&\quad +&\!\!\!\frac{2i}{m+2}tr_\theta d\ell(\cdot,J\cdot)\
\ -\frac{i}{m+2}tr_\theta d\ell(\cdot,J\cdot)\\[5mm]
&=\quad&\!\!\!\frac{2m+1}{m+1}\cdot scal^W\ +\ \frac{i}{m+2} tr_\theta d\ell(\cdot,J\cdot)\ \ \ .
\end{eqnarray*}
Thereby, we use again the torsion identities from section \ref{ab6}, which hold under the condition
of partial integrability.\hfill$\Box$\\

Theorem \ref{TH2} shows that the scalar curvature of the Fefferman metric $f_\theta$
and the Webster scalar curvature are proportional. For arbitrary $\ell\in\Omega^1(M;i\RR)$
this is not true.

\section{The Laplacian of the fundamental Killing vector} 
\label{ab8}

We discuss now properties of the fundamental vector field $S$ which is vertical along the 
$S^1$-fibre bundle $F_c$ in the Fefferman construction. It is easy to see that it is a Killing vector with
respect to any $\ell$-Fefferman metric $f_{\theta,\ell}$, i.e., $\eL_Sf_{\theta,\ell}=0$, where $\eL$ denotes
the Lie derivative.
However, the main goal of this section is to calculate the Laplacian of 
$S$. We are able to give an explicit expression for it.
The result is a first step in direction of a tractor calculus description for 
$\ell$-Fefferman spaces which come from partially integrable CR-spaces. We will complete 
the tractor description in the later sections. 

The fundamental vector field $S$ is uniquely determined by $A_{\theta,\ell}(S)=i\frac{m+2}{2}$
and it is lightlike with respect to the $\ell$-Fefferman metric $f_{\theta,\ell}$, where $\ell\in\Omega^1(M;i\RR)$.
Lemma \ref{LM6} shows that
\[f_{\theta,\ell}(\nabla^{f_{\theta,\ell}}_C S,B)=-f_{\theta,\ell}(\nabla^{f_{\theta,\ell}}_BS,C)\]
for all vectors $B,C\in TF_c$, i.e., $S$ is a Killing vector for any $\ell$-Fefferman metric $f_{\theta,\ell}$. 
Equivalently, for the dual $1$-form $\theta$ to $S$ on $(F_c,f_{\theta,\ell})$ holds 
\[\nabla^{f_{\theta,\ell}}\theta=1/2\cdot d\theta\ .\]
Now let $\Delta^{f_{\theta,\ell}}=d^*d+dd^*$ denote the Laplace-Beltrami operator acting on differential forms,
whereby $d^*$ is the codifferential with respect to $f_{\theta,\ell}$. 
By $tr_{f_{\theta,\ell}}\!\nabla^2$ we denote
the Bochner-Laplacian acting on arbitrary tensor fields $\rho$ through
\[tr_{f_{\theta,\ell}}\!\nabla^2\rho\ =\ 
\nabla^{f_{\theta,\ell}}_S\nabla^{f_{\theta,\ell}}_{T^*}\rho+\nabla^{f_{\theta,\ell}}_{T^*}\nabla^{f_{\theta,\ell}
}_{S}\rho
\ +\ \sum_{i=1}^{2m}\big(\ \nabla^{f_{\theta,\ell}}_{e_i^*}\nabla^{f_{\theta,\ell}}_{e_i^*}\rho
-\nabla^{f_{\theta,\ell}}_{\nabla^{f_{\theta,\ell}}_{e^*_i}e^*_i}\rho\ \big)\]
(with respect to our special choice of frame on $F_c$).
In general, for a Killing $1$-form 
$\theta$ holds \[d^*\theta=0\qquad\mbox{and}\qquad 
tr\nabla^2\theta=-1/2\cdot d^*d\theta=-1/2\cdot \Delta\theta\ .\]
By using the formulae of Lemma \ref{LM6} we find a simple expression for the Laplacian applied to
$\theta$ on $(F_c,f_{\theta,\ell})$.
\begin{PR} \label{PR1} (cf. \cite{Lei05}) Let $(F_c,f_{\theta,\ell})$ be the $\ell$-Fefferman space of a partially 
integrable CR-space 
$(M^n,T_{10})$  of dimension $n=2m+1$ with pseudo-Hermitian structure $\theta$ and $\ell\in\Omega^1(M;i\RR)$. Then 
\begin{enumerate} 
\item the fundamental 
$S$ is a Killing vector and the lift $\theta$ to $F_c$ of the pseudo-Hermitian form 
is the dual Killing $1$-form. In particular, 
\[\nabla^{f_{\theta,\ell}}\theta=1/2\cdot d\theta\ .\]
\item For the Laplace-Beltrami operator applied to $\theta$ on $(F_c,f_{\theta,\ell})$ holds
\[\Delta^{f_{\theta,\ell}}\theta=-2i\frac{(n-1)}{n+3}\cdot A_{\theta,\ell}+\left(\frac{scal^{f_{\theta,\ell}}}{n}
\ -\ \frac{2(n+1)}{n\cdot(n+3)}tr_\theta L_{d\ell}\right)\cdot\theta\]
resp. it holds
\[\mathcal{P}\theta=\frac{i}{n+3}\cdot A_{\theta,\ell}\ +\ \frac{(n+1)\cdot tr_\theta
L_{d\ell}}{n(n-1)\cdot(n+3)}\cdot\theta\ ,\]
where $\eP=\eP^{f_{\theta,\ell}}$ is the differential operator $\frac{1}{n-1}\left(\ 
tr_{f_{\theta,\ell}}\!\nabla^2+\frac{scal^{f_{\theta,\ell}}}{2n}\ \right)$. 
\end{enumerate}
\end{PR}
  
{\bf Proof.} Let $X,Y,Z$ denote sections in $H$ on $M$ such that their coordinates with respect to a local frame
$\{e_i:\ i=1,\ldots,2m\}$
are constant. It holds
\begin{eqnarray*}f_{\theta,\ell}(\nabla^{f_{\theta,\ell}}_{X^*}\nabla^{f_{\theta,\ell}}_{Y^*}S,S)
&=&-f_{\theta,\ell}(\nabla^{f_{\theta,\ell}}_{Y^*}S,
\nabla^{f_{\theta,\ell}}_{X^*}S)=-\frac{1}{4}f_{\theta,\ell} (\, (JY)^*,(JX)^*\, )\\[2mm]
&=&-\frac{1}{4}f_{\theta,\ell}(X^*,Y^*)\ ,\\[6mm]
f_{\theta,\ell}(\nabla^{f_{\theta,\ell}}_{X^*}\nabla^{f_{\theta,\ell}}_{Y^*}S,T^*)
&=&-f_{\theta,\ell}(\nabla^{f_{\theta,\ell}}_{Y^*}S,
\nabla^{f_{\theta,\ell}}_{X^*}T^*)\\[2mm]
&=&-\frac{1}{4}\big(\ L_\theta([X,T],JY)+L_\theta([JY,T],X)\ \big)\\
&&-i\frac{1}{2(m+2)}\Omega_{\theta,\ell}(JY,X),\\[6mm]
f_{\theta,\ell}(\nabla^{f_{\theta,\ell}}_{X^*}\nabla^{f_{\theta,\ell}}_{Y^*}S,Z^*)&=&
-f_{\theta,\ell}(\nabla^{f_{\theta,\ell}}_{Y^*}S,
\nabla^{f_{\theta,\ell}}_{X^*}Z^*)\\[1mm]&=&\frac{1}{2}\big(\ L_\theta(\nabla^W_XJZ,Y)+\eB_\theta(X,JZ,Y)\ \big).
\end{eqnarray*}
These formulae show that
\[\begin{array}{rcl} 
f_{\theta,\ell}(tr_{f_{\theta,\ell}}\nabla^2S,S)&=\quad\ &\!\!\!\!-\frac{m}{2}\\[6mm]
f_{\theta,\ell}(tr_{f_{\theta,\ell}}\nabla^2S,Z^*)&=\quad\ &\!\!\!\!\frac{1}{2}\sum_{i=1}^{2m} 
\big(\ L_\theta(\nabla^W_{e_i}JZ,e_i)+ L_\theta(\nabla^W_{e_i}e_i,JZ)\ 
\big)\\[2.5mm]
&\quad\ +&\!\!\!\!\frac{1}{2}\sum_{i=1}^{2m}\eB_\theta(e_i,JZ,e_i)\\[4mm]
&=\quad\ &\!\!\!\!\frac{1}{2}\sum_{i=1}^{2m}e_i\big(\ L_\theta(e_i,JZ)\ \big)\ =\ 0\\[7mm] 
f_{\theta,\ell}(tr_{f_{\theta,\ell}}\nabla^2S,T^*)&=\quad\ &\!\!\!\!-i\frac{1}{2(m+2)}\sum_{i=1}^{2m} 
\Omega_{\theta,\ell}(Je_i,e_i)\\[3.5mm]
&=\quad\ &\!\!\!\!-\frac{\ scal^W}{2(m+1)}+\frac{i}{2(m+2)}tr_\theta d\ell(\cdot,J\cdot)\\[3.5mm]
&=\quad\ &\!\!\!\!-\frac{1}{2n}scal^{f_{\theta,\ell}}+\frac{n+1}{n(n+3)}tr_\theta L_{d\ell}\ \ ,
\end{array}\]
whereby we use the relation 
\[scal^{f_{\theta,\ell}}=\frac{2m+1}{m+1}scal^W+\frac{1}{m+2}tr_\theta L_{d\ell}\ .\]
With the identity 
$f_{\theta,\ell}(tr_{f_{\theta,\ell}}\nabla^2S,\cdot)=tr_{f_{\theta,\ell}}\nabla^2\theta(\cdot)=-\frac{1}{2}
\Delta^{f_{\theta,\ell}}\theta (\cdot)$ 
we obtain the stated formula for the 
Laplace-Beltrami operator applied to $\theta$. \hfill$\Box$\\

\noindent
The importance of the next conclusion and the differential operator $\eP$ on a semi-Riemannian space 
in the context of conformal geometry will be explained in the next sections.

\begin{COR} Let $(F_c,f_{\theta,\ell})$ be the $\ell$-Fefferman space of $(M,T_{10})$
with respect to $\theta,\ell$ on $M$. It holds
\[\eP\theta=\frac{i}{n+3}A_{\theta,\ell}\qquad\mbox{if\ and only\ if}\qquad tr_\theta L_{d\ell}=0\ .\]
\end{COR}

Eventually, we define here the vector space
\[H_{tr}^1(M,T_{10}):=\{\ell\in\Omega^1(M;i\RR):\ tr_\theta d\ell(\cdot,J\cdot)=0\}/
\{\ell\in\Omega^1(M;i\RR):\ d\ell=0\} \]
to any partially integrable CR-space $(M,T_{10})$. 
This quotient can be understood as the vector space on which the affine space 
of those connection forms is modelled, which admit the same scalar curvature as the
Weyl connection on $F_c$ belonging to some fixed $\theta$, modulo gauge transformations. 
The class $c_\ell=[f_{\theta,\ell}]$ is locally conformally equivalent to any
$c_{\tilde{\ell}}=[f_{\theta,\tilde{\ell}}]$ with $\tilde{\ell}\in[\ell]$, where $[\ell]$ denotes 
the class of $1$-forms with purely imaginary values which differ from $\ell$ only by a closed form (i.e. gauge 
transformation).
Hence we can 
introduce the notion $c_{[\ell]}$ which denotes (locally) a uniquely determined conformal
structure, which we call the $[\ell]$-Fefferman conformal class on $F_c$.
The map
\[
\Psi:\ \ [\ell]\in H^1_{tr}(M,T_{10})\qquad\mapsto\qquad c_{[\ell]}=[f_{\theta,\ell}]
\]
is then bijective onto the space of local conformal structures on $F_c$
belonging to $\ell$-Fefferman metrics, which have the property
that
\[
\eP\theta=\frac{i}{n+3}A_{\theta,\tilde{\ell}}
\]
for all pseudo-Hermitian forms $\theta$ and $\tilde{\ell}\in[\ell]$.

\section{Aspects of conformal tractor calculus}%
\label{ab9}

We collect in this section notions and facts from
conformal tractor calculus. Thereby, we will restrict ourselves only 
to the very necessary parts of the apparatus, which we will need for our purposes in the remaining 
sections. 
We will omit the use of CR-tractors.
For a broader explanation of the topics that are raised
we refer to \cite{CG02}, \cite{Cap05a}, \cite{CS00}, \cite{CSS97}, \cite{CSS01}, \cite{Cap05c}.

Let $\frak{g}$ denote the Lie algebra $\frak{so}(2,n+1)$
belonging to the special orthogonal group $\SO(2,n+1)$, which acts by standard representation
on the Euclidean space $\RR^{2,n+1}$ with indefinite scalar product
\[\langle x,y\rangle=x_-y_++x_+y_-+(x_1,\ldots,x_{n+1})\ \mathbb{J}\ {}^t\!(y_1,\ldots,y_{n+1})\ ,\]
where the matrix \[\mathbb{J}=\left(\begin{array}{rc}-1& 0\\ 0 & I_{n}\end{array}\right)\]
gives rise to the Minkowski metric on $\RR^{n+1}$. (We choose here the special signature 
$(2,n+1)$ for $\frak{g}$, since
we want to work in the Lorentzian setting, which is related to strictly pseudoconvex CR-geometry.)
The Lie algebra $\frak{g}=\frak{so}(2,n+1)$ is $|1|$-graded:
\[\frak{g}=\frak{g}_{-1}\oplus\frak{g}_0\oplus \frak{g}_1\ ,\]
where $\frak{g}_0=\frak{co}(1,n)$, $\frak{g}_{-1}=\RR^{n+1}$ and $\frak{g}_1=\RR^{n+1*}$.
The $0$-part $\frak{g}_0$ decomposes further to its semisimple part $\frak{o}(1,n)$
and the center $\RR$.
We realise the subspaces $\frak{g}_0, \frak{g}_{-1}$ and $\frak{g}_{1}$ by matrices
of the form
\[\left(\begin{array}{ccc} 0& 0& 0\\m & 0& 0\\ 0& -{}^tm\mathbb{J}& 0
\end{array}\right)\in\frak{g}_{-1}\ ,\quad\ \
\left(\begin{array}{ccc} -a& 0& 0\\0 & A& 0\\ 0& 0& a
\end{array}\right)\in\frak{g}_0\ ,\quad\ \
\left(\begin{array}{ccc} 0& l& 0\\0 & 0& -\mathbb{J}\ {}^tl\\ 0& 0& 0
\end{array}\right)\in\frak{g}_{1}\ .\]
The commutators with respect to these matrices are given by
\[ \begin{array}{ll}
{[\ ,\ ]}:\ \frak{g}_0\times\frak{g}_0\to \frak{g}_0\ ,\qquad &
[(A,a),(A',a')]=(AA'-A'A,0)\\[2mm]
{[\ ,\ ]}:\ \frak{g}_0\times \frak{g}_{-1}\to \frak{g}_{-1}\ ,\qquad &   
[(A,a),m]=Am+am\\[2mm]
{[\ ,\ ]}:\ \frak{g}_1\times \frak{g}_{0}\to \frak{g}_{1}\ ,\qquad &
[l,(A,a)]=lA+al\\[2mm]
{[\ ,\ ]}:\ \frak{g}_{-1}\times \frak{g}_{1}\to \frak{g}_{0}\ ,\qquad &
[m,l]=(ml-\mathbb{J}\ {}^t(ml)\mathbb{J},lm)\ ,\end{array}\qquad \]
where $(A,a)$, $(A',a')\in \frak{o}(1,n)\oplus\RR,\ m\in\RR^{n+1},\ l\in
\RR^{n+1*}$.

The subalgebra
\[\frak{p}:=\frak{g}_0\oplus\frak{g}_1\]
of $\frak{g}$ is a parabolic, i.e., it contains a maximal solvable subalgebra of
$\frak{g}$. We also denote $\frak{p}_+:=\frak{g}_{1}$. 
While the grading of $\frak{g}$ is not $\frak{p}$-invariant,
it gives rise to an invariant filtration
\[\frak{g}\supset\frak{p}\supset\frak{p}_+\ .\]
The subgroup $P$ of $G=\SO(2,n+2)$, whose adjoint action on $\frak{g}$ preserves
this filtration is a parabolic subgroup with Lie algebra $\frak{p}$. And the subgroup
$G_0$ which preserves the grading of $\frak{g}$ is the group $\CSO(1,n)=\SO(1,n)\times\RR^+$ with Lie algebra
$\frak{g}_0=\frak{co}(1,n)$. Moreover, the exponential map restricts to a diffeomorphism from $\frak{p}_+$
onto a normal subgroup $P_+\subset P$ and the parabolic $P$ is the semidirect product of $G_0$ and $P_+$.
Eventually, we note that the homogeneous space $G/P$ is the flat model of conformal geometry in Lorentzian signature
$(1,n+1)$.

The Killing form of
$\frak{g}$ defines a duality between $\frak{g}/\frak{p}$ and $\frak{p}_+$
as $\frak{p}$-modules (with respect to the adjoint action). In particular, for each $k\in\mathbb{N}$ 
we get a $\frak{p}$-module isomorphism 
\[H\!om(\Lambda^k\frak{g}/\frak{p},\frak{g})\cong \Lambda^{k}\frak{p}_+\otimes \frak{g}\ .\]
The latter spaces are the groups in the 
standard complex computing the Lie algebra homology of $\frak{p}_+$ with coefficients in
$\frak{g}$. The differentials in this standard complex are linear maps denoted by
\[\partial^*:H\!om(\Lambda^k\frak{g}/\frak{p},\frak{g})\to H\!om(\Lambda^{k-1}\frak{g}/\frak{p},\frak{g})\ .\]
For $k=1$ the explicit formula for $\partial^*$ is given by
\[\phi\in\ \frak{p}_+\otimes\frak{g}\qquad\mapsto\qquad 
\partial^*\phi=\sum_{i=1}^{n+1}[\eta_i,\phi(\xi_i)]\ \in\frak{g}\ \ ,\]  
where $\{\xi_i:\ i=1,\ldots,n+1\}$ is some basis of $\frak{g}_{-1}\cong\frak{g}/\frak{p}$ 
and $\{\eta_i:\ i=1,\ldots,n+1\}$ is the dual basis of $\frak{p}_+$.
By construction, it is $\partial^*\circ\partial^*=0$ for all $k$ and we denote the $k$th homology group 
$ker(\partial^*)/Im(\partial^*)$ which is a 
$\frak{p}$-module by $\mathcal{H}^k_\frak{g}$. It is well known that $\frak{p}_+$ acts trivially
on $\mathcal{H}^k_\frak{g}$ (cf. \cite{Cap05a}).   

Now let $(F^{n+1},c)$ be a smooth orientable manifold  with conformal structure $c$ of signature $(1,n)$. 
The conformal
structure is given by a $G_0$-principal bundle reduction $\mathcal{G}_0(F)$
of the general linear frame bundle over $F$. By the process of prolongation we obtain
from $\mathcal{G}_0(F)$ the $P$-principal fibre bundle $\mathcal{P}(F)$ which is a reduction
of the second order linear frame bundle on $F$ and again determines the conformal structure $c$ on $F$
uniquely (cf. e.g. \cite{Kob72} \cite{CSS97}). A Cartan connection $\omega$ is a smooth $1$-form on $\mathcal{P}(F)$ 
with 
values
in $\frak{g}$ such that
\begin{enumerate}
\item $\omega(\chi_A)=A\ \ $ for all fundamental fields $\chi_A$, $A\in\frak{p}$,
\item $R_g^*\omega=Ad(g^{-1})\circ\omega\ \ $ for all $g\in P$\ \ and
\item $\omega|_{T_u\mathcal{P}(F)}:T_u\mathcal{P}(F)\to\frak{g}\ \ $ is a linear isomorphism for all
$u\in\mathcal{P}(F)$.
\end{enumerate}

The curvature $2$-form $\Omega$ of a Cartan connection $\omega$ is given by
\[\Omega=d\omega+\frac{1}{2}[\omega,\omega]\ .\]
The corresponding curvature function $\kappa:\mathcal{P}(F)\to H\!om(\Lambda^2\frak{g},\frak{g})$ 
is given by
\[\kappa(u)(X,Y):=d\omega(\omega_u^{-1}(X),\omega_u^{-1}(Y))+[X,Y]\]
for $u\in\mathcal{P}(F)$ and $X,Y\in\frak{g}$. The curvature function $\kappa$ vanishes if one of its entries lies in 
$\frak{p}\subset\frak{g}$. Hence we can view
$\kappa$ as an $P$-equivariant smooth function \[\kappa:\mathcal{P}(F)\to H\!om(\Lambda^2\frak{g}/\frak{p},\frak{g})\ .\]
An important fact of conformal geometry is that there exists a distinguished
Cartan connection $\omega_{nor}$ on $\mathcal{P}(M)$, which we call the normal Cartan connection of conformal
geometry (cf. e.g. \cite{Kob72}, \cite{CSS97}). The normal Cartan connection $\omega_{nor}$ is uniquely determined by 
the condition
\[\partial^*\circ\kappa=0\ .\]     

Furthermore, we can extend the $P$-principal fibre bundle $\mathcal{P}(F)$ by the structure group $G$. This gives
rise to the $G$-principal fibre bundle $\mathcal{G}(F)=\mathcal{P}(F)\times_{P}G$ over $F$. The normal Cartan connection
$\omega_{nor}$ on $\mathcal{P}(F)$ extends by equivariance as well and becomes a principal fibre bundle connection on 
$\mathcal{G}(F)$ which we denote again by $\omega_{nor}$. Now let $\mathbb{V}$ be a $G$-representation.
From the representation $\mathbb{V}$ we obtain 
the vector bundle $\mathcal{V}:=\mathcal{G}(F)\times_{G}\mathbb{V}$. We call $\mathcal{V}$ the tractor bundle
which belongs to the $G$-representation $\mathbb{V}$.
The normal connection $\omega_{nor}$ on $\mathcal{G}(F)$ induces on each tractor bundle $\mathcal{V}$
an invariant covariant derivative:
\[\nabla^{nor}:\Gamma(TF^*)\otimes\Gamma(\mathcal{V})\to\Gamma(\mathcal{V})\ .\]

Of particular interest for our purpose are the standard and adjoint tractor bundles of a space
$(F,c)$ with conformal structure. We denote by $\mathcal{T}(F):=\mathcal{G}(F)\times_{G}\RR^{2,n+1}$  
the standard tractor bundle and by  $\mathcal{A}(F):=\mathcal{G}(F)\times_{Ad(G)}\frak{g}$ the adjoint tractor
bundle. In the latter case the adjoint representation of $G$ on $\frak{g}$ can be restricted to $P$,
which preserves the filtration of $\frak{g}$. Hence we obtain a natural filtration
\[\mathcal{A}(F)\supset\mathcal{A}^0(F)\supset\mathcal{A}^{1}(F)\ ,\]
whereby $\mathcal{A}^0(F)=\mathcal{P}(F)\times_{Ad(P)}\frak{p}$ and  the 
bundle $\mathcal{A}^1(F)=\mathcal{P}(F)\times_{P}\frak{p}_+$
is the dual $T^*F$ of the tangent bundle of $F$. Notice that with this
identification the dual tangent bundle $T^*F$ is canonically contained in $\mathcal{A}(F)$. 
Furthermore, by choosing a Lorentzian metric $g$ in $c$ we can reduce
the bundles $\mathcal{G}(F)$ and $\mathcal{P}(F)$
(by using the corresponding Weyl structure of $g$)
to the structure group 
$\SO(1,n)$, which is the semisimple part of $G_0$. This reduction 
gives rise
to a grading on the adjoint tractor bundle. It is
\[\mathcal{A}(F)\ \cong_g\ TF\oplus\frak{co}(TF)\oplus T^*F\ .\]
Thereby,  notice that $TF=\mathcal{G}_0(F)\times_{G_0}\frak{g}_{-1}$ 
and $\frak{co}(TF)$ is the subset of the endomorphism bundle 
$End(TF)$ consisting of skew-symmetric maps plus non-zero multiples of the identity map.
The concrete identification of $\mathcal{A}(F)$ with this graded sum depends 
very much on the choice of  $g$ in $c$. The standard tractor bundle splits as graded sum
with respect to a metric $g$ in $c$ by
\[\mathcal{T}(F)\cong_g\RR\oplus TF\oplus \RR\ .\]

The Lie bracket on $\frak{g}$ is $Ad(P)$-invariant. Hence there is an induced algebraic bracket
on the adjoint tractor bundle $\mathcal{A}(F)$, which we denote by 
$\{\cdot,\cdot\}$. This bracket admits with respect to some metric $g$ in the 
conformal class following expressions. For elements $(\xi,\varphi,\omega), 
(\tau,\psi,\eta)\in\mathcal{A}_p(F)$ at $p\in F$,
whereby the components are given with respect to the induced grading of $\mathcal{A}(F)$ by $g$, it holds
\[\begin{array}{c}
\{\ (\xi,\varphi,\omega)\ ,\ (\tau,\psi,\eta)\ \}\ =\\[1mm] 
(\ \{\xi,\psi\}+\{\varphi,\tau\}\ ,\ \{\xi,\eta\}+\{\varphi,\psi\}+\{\omega,\tau\}\ ,\ 
\{\varphi,\eta\}+\{\omega,\psi\}\ )\ .\end{array}
\]
Moreover, it is
\[\begin{array}{rcl}
\{e_i,\psi\}&=&-\ \psi(e_i)\\[2mm]
\{e_i^*,\psi\}&=&\quad e_i^*\circ\psi\\[2mm]
\{e_i,e_j^*\}&=&\quad e_i^*\otimes e_j^*\ -\ \varepsilon_i\varepsilon_je_j\otimes e_i^*\qquad\mbox{for}\ \ i\neq j 
\\[2mm]
\{e_i,e_i^*\}&=&\quad id|_{TF}\\[2mm]
\{e_i,\eta\}&=&\quad e_i\otimes\eta\ -\ \sum_{j=1}^{n+1}\varepsilon_i\varepsilon_j\eta(e_j)e_j\otimes 
e_i^*\ +\ \eta(e_i)id\\[2mm]
\{\xi,e_j^*\}&=&\quad \xi\otimes e_j^*\ -\ \varepsilon_j e_j\otimes g(\xi,\cdot)\ +\ e_j^*(\xi)id\ ,
\end{array}\]
where $\{e_i:\ i=1,\ldots,n+1\}$ denotes an orthonormal frame on $F$ with respect to $g$ and 
$\varepsilon_i:=g(e_i,e_i)$.  

Furthermore, the $\frak{g}$-action on $\RR^{2,n+1}$ which is compatible with the $P$-action
gives rise to a multiplication on $\mathcal{T}(F)$ by adjoint tractors, which we denote by
\[\begin{array}{cccl}
\bullet:& \mathcal{A}(F)\otimes\mathcal{T}(F)&\to& \mathcal{T}(F)\quad .\\[2mm]
&(A,t)&\mapsto& A\bullet t 
\end{array}
\]
In this respect, we can view an adjoint tractor as endomorphism on standard tractors. We also use 
the notation $A^2\in End(\mathcal{T}(F))$ for the superposition $A\bullet A$ of $A\in\mathcal{A}(F)$
acting on standard 
tractors.
Of course, $A^2$ is in general not an adjoint tractor.

Since the differential 
$\partial^*:\frak{p}_+\otimes\frak{g}\to\frak{g}$
is $P$-equivariant, it induces a vector bundle map 
\[\partial^*:T^*F\otimes\mathcal{A}(F)\to \mathcal{A}(F)\ .\]
The image of $\partial^*$ is $\frak{p}$ resp. $\mathcal{A}^0(F)$. Let us denote 
$H^1_\frak{g}F:=\mathcal{P}(F)\times_P\mathcal{H}_\frak{g}^1$.
It holds 
\[H^1_\frak{g}F=\mathcal{A}(F)/\mathcal{A}^0(F)\ .\]
We mentioned above that the group $P_+$ acts trivially on $\mathcal{H}_\frak{g}^1$. For that reason we can
view $H_\frak{g}^1F$ as the associated bundle $\mathcal{G}_0\times_{G_0}\mathcal{H}_\frak{g}^1$, which 
is naturally identified with the tangent bundle of $F$:
\[H^1_\frak{g}F\ \cong\ TF\ .\] 
Let $\pi_H:\mathcal{A}(F)\to H^1_\frak{g}(F)$ denote the natural projection. This projection does not depend on
the choice of a metric or a Weyl structure in $c$. If $\alpha$ is a section of $\mathcal{A}(F)$ then
we can consider  
$\pi_H(\alpha)$ as a tangent vector field, which we often call the 'first slot' of the adjoint tractor $\alpha$. 

On any space $(F,c)$ with conformal structure
there exists a certain (conformally) invariant differential operator, which assigns
to any section
of the bundle $H_{\mathbb{V}}^kF$ of homology groups of $\frak{p}_+$ with coefficients in some
$G$-module $\mathbb{V}$ an element of $\Omega^k(F,\mathcal{V})$, where the latter  
is the space of smooth
$k$-forms with
values in the tractor bundle belonging to the representation $\mathbb{V}$. This operator is called
the splitting operator and plays a prominent role in the construction of BGG sequences (cf. 
\cite{CSS01}, \cite{Cap05a}).
We denote this operator by $\mathcal{S}$. For the rest of the section, we aim to calculate 
$\eS$ explicitly in terms of a
metric $g$ in the conformal class $c$ on $F$
for the case when it is a map from vector fields
to sections of the adjoint tractor bundle:
\[\mathcal{S}: \frak{X}(F)\to \Gamma(\mathcal{A}(F))\ .\]

In general, the splitting operator is determined by the following property (cf. \cite{CSS01}).
Let $d^{nor}:\Omega^k(F,\mathcal{V})\to \Omega^{k+1}(F,\mathcal{V})$ denote the covariant exterior derivative,
which is induced by the normal connection on the $G$-principal fibre bundle $\mathcal{G}(F)$.
Then for arbitrary
$\tau\in\Gamma(H^k_{\mathbb{V}}F)$, it holds
\begin{enumerate}
\item $\partial^*(\mathcal{S}(\tau))=0$\quad  and\quad $\pi_H(\mathcal{S}(\tau))=\tau$,\\[-2mm]
\item $\partial^*(d^{nor}\mathcal{S}(\tau))=0$ .
\end{enumerate}
In the particular  situation when $k=0$ and $\mathbb{V}=\frak{g}$ the covariant exterior derivative is
just equal to $\nabla^{nor}$. Hence the splitting operator
\[\mathcal{S}:\frak{X}(F)\to \Gamma(\mathcal{A}(F))\]
is determined by the property that for all vector fields $\tau$ on $F$, it holds
\[\pi_H(\mathcal{S}(\tau))=\tau\qquad\mbox{and}\qquad \partial^*(\nabla^{nor}\mathcal{S}(\tau))=0\ .\]
Moreover, in this situation the map $\mathcal{S}$ admits the property that if $\tau\in\frak{X}(F)$ is
a conformal Killing vector field on $(F,c)$ then
\[ \nabla^{nor}_X\mathcal{S}(\tau)=-\Omega^{nor}(\tau,X)\qquad\mbox{for\ all\ }X\in TF\  .\]
Conversely, it is also known that the 'first slot' of any adjoint tractor, which solves the latter
equation, is a conformal Killing vector field (cf. \cite{CD01}, \cite{Cap05c}). 

The covariant exterior derivative $d^{nor}=\nabla^{nor}$
acting on a section $(\tau,\psi,\eta)$ in 
$\mathcal{A}(F)$
with respect to a metric $g$ in the conformal class $c$
on $F$ is given by 
\begin{eqnarray*}
d^{nor}(\tau,\psi,\eta)\!\!&=&\!\!(\ \nabla^g_\cdot\tau\ ,\ \nabla^g_\cdot\psi\ ,\
\nabla^g_\cdot\eta\ )+\{\ (\cdot,0,P^g(\cdot))\ ,\ (\tau,\psi,\eta)\ \}\\
&=&\!\!(\ \nabla^g_\cdot\tau\ ,\ \nabla^g_\cdot\psi\ ,\ \nabla^g_\cdot\eta\ )+(\{\cdot,\psi\}\ ,\ 
\{\cdot,\eta\}+\{P^g(\cdot),\tau\}\ ,\
\{P^g(\cdot),\psi\}),
\end{eqnarray*}
where $\nabla^g$ denotes the Levi-Civita connection and
\[\begin{array}{cccl}P^g:& TF&\to& T^*F\\[1mm]
&X&\mapsto&\frac{1}{n-1}\left(\ \frac{scal^g}{2n}g(X,\cdot)-Ric^g(X,\cdot)\ \right)\end{array}
\]
is the Schouten tensor with respect to $g$ (cf. \cite{CG02}). Applying the differential $\partial^*$ results to
\[\partial^*d(\tau,\psi,\eta)=\left(\begin{array}{c}0\\[2mm]
\sum_{i=1}^{n+1}(\ \{e_i^*,\nabla_{e_i}^g\tau\}+\{e_i^*,\{e_i,\psi\}\}\ )\\[2mm]
\sum_{i=1}^{n+1}(\ 
\{e_i^*,\nabla^g_{e_i}\psi\}\ +\ \{e_i^*,\{e_i,\eta\}\}\ +\ \{e_i^*,\{P^g(e_i),\tau\}\}\ )\end{array}\right)\ .
\]
With the formulae for the bracket $\{\cdot,\cdot\}$  from above we calculate
\[\begin{array}{rcl}
\sum_{i=1}^{n+1}\{e_i^*,\{e_i,\psi\}\}&=&\quad\psi\ +\ tr_g\psi\cdot id\ -\ \sum_{i=1}^{n+1}\varepsilon_i 
g(\psi(e_i),\cdot)\\[2mm]
\sum_{i=1}^{n+1}\{e_i^*,\nabla^g_{e_i}\tau\}&=&-\ \nabla^g_\cdot\tau\ -\ 
tr_g\nabla^g\tau\ +\ \sum_{i=1}^{n+1}\varepsilon_i
g(\nabla^g_{e_i}\tau,\cdot)\\[2mm]
\sum_{i=1}^{n+1}\{e_i^*,\{e_i,\eta\}\}&=&\quad(n+1)\cdot \eta\\[2mm]
\sum_{i=1}^{n+1}\{e_i^*,\nabla^g_{e_i}\psi\}&=&\quad\sum_{i=1}^{n+1}e_i^*\circ\nabla^g_{e_i}\psi\\[2mm]
\sum_{i=1}^{n+1}\{e_i^*,\{P^g(e_i),\tau\}\}&=&-\ 2P^g(\tau)\ +\ tr_gP^g\cdot g(\tau,\cdot),
\end{array}
\]
where $tr_gP^g=-\frac{scal^g}{2n(n+1)}$.
One can immediately see that the equation $\partial^*d(\tau,\psi,\eta)=0$ is solved for an arbitrary
vector field $\tau$ on $F$ by setting
\[\begin{array}{rcl}
\psi&=&\quad \nabla^g\tau\qquad\quad\mbox{and}\\[2mm]
\eta&=&-\ \frac{1}{n+1}\left(\ (\sum_{i=1}^{n+1}e_i^*\circ\nabla^g_{e_i}\psi)-2P^g(\tau)+tr_gP^g\cdot g(\tau,\cdot)\
\right)\ .
\end{array}\]
By inserting the first equation into the second one we obtain
\[\begin{array}{rcl}
\eta&=&-\frac{1}{n+1}\left(\ (\sum_{i=1}^{n+1}e_i^*\circ\nabla^g_{e_i}(\nabla^g\tau)(\cdot))
-2P^g(\tau)+tr_gP^g\cdot g(\tau,\cdot)\ \right)\ .
\end{array}\]

We reformulate the latter  expression for $\eta$ in the case when $\tau$ is a conformal Killing vector field.
In this situation, it holds 
\begin{eqnarray*}
\sum_{i=1}^{n+1}e_i^*\circ\nabla^g_{e_i}(\nabla^g\tau)&=&-g(tr_g\nabla^2\tau,\cdot)+\frac{2}{n+1}d(div(\tau))
\qquad\qquad\mbox{and}\\
d(div(\tau))&=&-g(tr_g\nabla^2\tau,\cdot)+\frac{2}{n+1}d(div(\tau))-Ric^g(\tau,\cdot)\ ,
\end{eqnarray*}
where $div(\tau)=tr_gg(\nabla^g_\cdot\tau,\cdot)$.
This implies 
\begin{eqnarray*}
\frac{n-1}{n+1}\cdot d(div(\tau))&=&-g(tr_g\nabla^2\tau,\cdot)\ -\ Ric^g(\tau,\cdot)\ ,\\
\sum_{i=1}^{n+1}e_i^*\circ\nabla^g_{e_i}(\nabla^g\tau)&=&-\frac{n+1}{n-1}g(tr_g\nabla^2\tau,\cdot)\ -\
\frac{2}{n-1}Ric(\tau,\cdot)\ ,
\end{eqnarray*}
and finally, we obtain
\[\eta=\frac{1}{n-1}\cdot g(tr_g\nabla^2\tau,\cdot)\ +\ \frac{scal^g}{2n\cdot (n-1)}g(\tau,\cdot)\ .\]
We conclude that the splitting operator $\mathcal{S}$ applied to a conformal Killing vector $\tau$ on $(F,c)$
is given with respect to some metric $g$ in $c$ by
\[
\mathcal{S}(\tau)=(\ \tau\ ,\ \nabla_\cdot^g\tau\ ,\ \mathcal{P}(g(\tau,\cdot))\ )\ \ ,
\]
where the differential operator $\eP=\eP^g$ is defined as in Proposition \ref{PR1} by
\[\mathcal{P}:=\frac{1}{n-1}\left(\
tr_g\nabla^2+\frac{scal^g}{2n}\ \right)\quad .\]

\section{Conformal tractor calculus for Fefferman space}
\label{ab10}

Now we come back to the realm of Fefferman spaces and
study their conformal tractor calculus.
It is well known in case of integrable
CR-structures that the standard tractor bundle of the
Fefferman conformal class admits a complex structure $J_{CR}$, which
arises naturally as the lift of the invariant complex structure
on the standard CR-tractors. 
However, the complex structure $J_{CR}$ has another important property,
namely it is parallel with respect to the normal conformal tractor connection $\nabla^{nor}$ and
the corresponding conformal vector field to $J_{CR}$ (i.e., the 'first slot' of $J_{CR}$) is exactly twice the
fundamental Killing field $S$ in fibre direction on $F_c$. 
We show here that a similar statement is true in the partially integrable case, but for general $\ell$-Fefferman
spaces the situation changes. 

Let $(F^{n+1}_c,f_{\theta,\ell})$ be a $\ell$-Fefferman space over a strictly pseudoconvex, 
partially integrable  CR-space $(M^n,T_{10})$ with pseudo-Hermitian structure $\theta$
and $\ell\in\Omega^1(M;i\RR)$. 
Let $\mathcal{T}(F_c)$
be the standard tractor bundle to the conformal structure $[f_{\theta,\ell}]$ on $F_c$.
Moreover, let $\mathcal{A}(F_c)$ be the adjoint tractor bundle,
which splits to a graded sum by the choice of $f_{\theta,\ell}$. Let us denote $R:=2S$.
We define the section $J_{CR}$ in  $\mathcal{A}(F_c)$ with respect to the graded sum corresponding to 
$f_{\theta,\ell}$ as
\[J_{CR}\ :=\ (\ R\ ,\ J_{\theta,\ell}\ ,\ \frac{2i}{n+3}A_{\theta,\ell}\ )\ ,\] 
whereby $J_{\theta,\ell}$ denotes the horizontal lift 
of the almost complex structure $J$ on $H$ 
to $TF_c$ (with respect to $A_{\theta,\ell}$). 
\begin{PR} \label{PR2} (1) The endomorphism \[J_{CR}=(R,J_{\theta,\ell},\frac{2i}{n+3}A_{\theta,\ell})\] 
defined as section  
of the adjoint tractor bundle $\mathcal{A}(F_c)$ does not depend on the choice
of $\theta$, i.e., $J_{CR}$ is a CR-invariant.\\[1mm]
(2) It is \[J^2_{CR}=-id|_{\mathcal{T}(F_c)}\ ,\] i.e., $J_{CR}$ is a complex structure on the standard tractor
bundle $\mathcal{T}(F_c)$. 
\end{PR}
{\bf Proof.} (1) In fact, the complex structure $J_{CR}$ is the lift of the invariant complex structure on
the CR-standard tractors of $M$. Since we do not introduce CR-tractors, we can not use this argument here.
Alternatively, we could apply the transformation law for the graded sum of adjoint tractors under rescaling of the 
metric (cf. \cite{CG02}) to show that the definition of $J_{CR}$ does not depend on $\theta$. However, we do not aim 
to 
introduce this transformation law either. So we postpone the proof of invariance
until Proposition \ref{PR3}. There, we will see that $J_{CR}$ is even a conformally invariant object
on $\mathcal{A}(F_c)$.\\
(2) With respect to the metric $f_{\theta,\ell}$ we can describe the action of $J_{CR}$ on a
standard tractor $t=(a,\xi,b)$ with $a\in \RR$, $\xi\in TF_c$ and $b\in\RR$ as application of the matrix
\[J_{CR}=
\left(
\begin{array}{ccc} 
\mbox{}\quad 0\qquad& \frac{2i}{n+3}A_{\theta,\ell}& 0\\[2mm]
\mbox{}\quad R\qquad& J_{\theta,\ell}& \frac{-2i}{n+3}A^b_{\theta,\ell}\\[2mm]
\mbox{}\quad 0\qquad&-f_{\theta,\ell}(R,\cdot)&0
\end{array}
\right)\quad ,
\]
where $A^b_{\theta,\ell}$ denotes the dual vector to $A_{\theta,\ell}$ with respect 
to $f_{\theta,\ell}$. For the middle entry of this matrix note that 
$f_{\theta,\ell}(\nabla_X^{f_{\theta,\ell}}R,Y)=d\theta(J_{\theta,\ell}X,Y)$ for all $X,Y\in TF_c$. 
Then it is \[J_{CR}\bullet t=
\left(\begin{array}{c}
\frac{2i}{n+3}A_{\theta,\ell}(\xi)\\[2mm]
a\cdot R+J_{\theta,\ell}(\xi)-\frac{2ib}{n+3}A^b_{\theta,\ell}\\[2mm]
-f_{\theta,\ell}(R,\xi)\end{array}
\right)\quad ,
\]
where $J_{\theta,\ell}$ acts trivially on the complement of the horizontally lifted distribution $H$.
Applying $J_{CR}$ by multiplication again results to
\[
J_{CR}\bullet J_{CR}\bullet t=
\left(
\begin{array}{c}
\frac{2ia}{n+3}A_{\theta,\ell}(R)\\[2mm]
\frac{2i}{n+3}A_{\theta,\ell}(\xi)\cdot 
R-pr_H(\xi)+\frac{2i}{n+3}f_{\theta,\ell}(R,\xi)A^b_{\theta,\ell}\\[2mm]
\frac{2ib}{n+3}A_{\theta,\ell}(R)
\end{array}
\right)\ =\ -t\quad ,
\]
whereby we use  $\frac{2i}{n+3}A_{\theta,\ell}(R)=-1$ and $pr_H$ denotes the projection of $TF_c$ to 
the horizontal lift of $H$.
\hfill$\Box$\\

We know so far that $J_{CR}$ is a complex structure on $\eT(F_c)$, whose 'first slot' 
is the (conformal) Killing vector field $R$, which is vertical in the fibres of the $\ell$-Fefferman construction.
With the help of the splitting operator $\eS$ we obtain the conformal invariant $\eS(R)$ on 
the adjoint tractor bundle $\mathcal{A}(F_c)$. We also set \[\mathcal{U}:=\frac{n+1}{n\cdot(n-1)(n+3)}
(0,0,tr_\theta L_{d\ell}\cdot\theta)\in\mathcal{A}(F_c)\ .\] Note that the $1$-form
$tr_\theta L_{d\ell}\cdot\theta$ does not depend on the pseudo-Hermitian structure $\theta$, i.e.,
$tr_\theta L_{d\ell}\cdot\theta$ is a uniquely determined section in $T^*F_c$. Since
$T^*F_c$ is canonically included in the adjoint tractor bundle, we know that $\mathcal{U}$ is a
uniquely defined section in $\mathcal{A}(F_c)$ and does not depend on the choice of $f_{\theta,\ell}$.
\begin{PR} \label{PR3} Let $(F^{n+1}_c,[f_{\theta,\ell}])$ be the $\ell$-Fefferman space of a 
strictly pseudoconvex, partially integrable CR-space
$(M^n,T_{10})$ with pseudo-Hermitian form $\theta$, correction term $\ell\in\Omega^1(M;i\RR)$ 
and fundamental Killing 
vector field $R$. It holds
\begin{enumerate}
\item
\[\mathcal{S}(R)=J_{CR}\ +\ \mathcal{U}\quad .\]
In particular, $J_{CR}$ is a conformal invariant on $\mathcal{A}(F_c)$. 
\item
It holds
\[ \eS(R)=J_{CR}\qquad\mbox{if\ and\ only if}\qquad [\ell]\in H^1_{tr}(M,T_{10})\quad . \]
In this case it is\[\nabla^{nor}_{\cdot}J_{CR}=-\Omega^{nor}(R,\ \cdot\ )\ \ .\]
\end{enumerate}
\end{PR}
Proposition \ref{PR3} follows immediately from Proposition \ref{PR1} and the definition of $J_{CR}$ and 
$\mathcal{U}$. We note that the curvature expression $\Omega^{nor}(R, \cdot )$ is entirely determined
by the Nijenhuis tensor $\eN$ and the $1$-form $\ell$.
The result also shows the independence of $J_{CR}$ from the choice of a pseudo-Hermitian form. This
proves the first statement of Proposition \ref{PR2}.

\section{The reconstruction}
\label{ab11}

In Proposition \ref{PR3} of the previous section we have seen how to obtain through the $\ell$-Fefferman construction
conformal structures which admit an orthogonal complex structure $\eR$ on the standard tractor bundle whose 
'first slot' $R$ is a conformal vector field. We want to argue now in the reversed 
direction
and discuss a reconstruction result from the existence of such a complex structure $\eR$. Our discussion
will be of local nature, since we do not want to assume a $S^1$-fibration on our initial space.

Let $(F^{n+1},c)$ be a smooth manifold of dimension $n+1$ with conformal structure $c$.
We assume here that there exists an adjoint tractor $\eR\in\Gamma(\mathcal{A}(F))$ which acts as
complex structure on the standard tractors $\mathcal{T}(F)$ such that the
'first slot' $\pi_H(\eR)=R$ is a conformal vector field on $F$, i.e.,
it holds
\[
\eR^2=-id|_{\eT(F)}\qquad\mbox{and}\qquad\nabla^{nor}_\cdot\eR=-\Omega^{nor}(R,\cdot)\ .
\]
The splitting operator applied to $R$ reproduces the adjoint tractor: $\eS(R)=\eR$.
In the following discussion, we will show
step by step how to construct from the existence of $\eR$ a (local) CR-structure on a quotient manifold.
The following Lemma \ref{LM8} is the first step.

\begin{LM} \label{LM8}
Let 
\[
\beta=\left(\
\begin{array}{ccc}
-a& l&0\\[2mm]
m& A& -\mathbb{J}\ {}^tl\\[2mm]
0&-{}^tm\mathbb{J}& a
\end{array}\right)
\]
be an arbitrary matrix in $\frak{g}=\frak{so}(2,n+1)$ with $\beta^2=-id$.
Then the vectors $m$ and $\mathbb{J}\ {}^tl\in\RR^{n+1}$ are non-zero and lightlike.
Moreover, it is $A(m)=am$ and $l\circ A=al$ and the restriction of $A$ to the $\mathbb{J}$-orthogonal complement 
$W$ of
$span\{ m,\mathbb{J}{}^tl \}$ in $\RR^{n+1}$ is a complex structure, i.e., $A^2|_W=-id$.\\
\end{LM}
{\bf Proof.} It is
\[
\beta^2=\left(\
\begin{array}{ccc}
a^2+lm& -al+lA&-l\mathbb{J}\ {}^tl\\[2mm]
-am+Am&ml+A^2+\mathbb{J}\ {}^tl{}^tm\mathbb{J}& -A\mathbb{J}\ {}^tl-a\mathbb{J}\ {}^tl\\[2mm]
-{}^tm\mathbb{J}m&-{}^tm\mathbb{J}A-a{}^tm\mathbb{J}& {}^tm{}^tl+a^2
\end{array}\right)\quad .
\]
This matrix square immediately shows that the condition $\beta^2=-id$ implies the statements of Lemma $\ref{LM8}$.
\hfill$\Box$\\

In general, with respect to a metric $f$ in $c$ and a corresponding orthonormal frame 
the matrix $\beta\in\frak{g}$ represents an adjoint tractor at some point $p\in F$. In particular, 
since we have chosen a metric $f$, the vector $m$ which determines the $\frak{g}_{-1}$-part 
of the matrix $\beta$ represents the 'first slot' of that adjoint tractor. Moreover, the matrix $A$
gives rise to a skew-symmetric endomorphism on $T_pF$ and $l$ is a lightlike $1$-form at $p\in F$.   
For our particular adjoint tractor $\eR$ this implies the following objects. It is $\eR=(\ R\ ,\ \nabla^fR\ ,\ 
\eP(f(R,\cdot))\ )$,
where the 'first slot' $R$ is a non-zero, lightlike conformal Killing vector field and the $1$-form $\eP(f(R,\cdot))$ 
is 
non-zero and lightlike as well. The skew-symmetric part of the endomorphism $\nabla^fR$ restricts on the 
$f$-orthogonal complement
$W^f$ of the lightlike vectors $R$ and the dual of $\eP(f(R,\cdot))$ in $TF$ to an
orthogonal complex structure $J^f$. 

Now we denote by
$\theta^f$ the dual of $R$ with respect to $f$ in $c$. The subbundles $\RR R$ and $E:=ker(\theta^f)$
give rise to a flag in $TF$. The quotient $Q:=E/\RR R$ is via $f$ identified with 
the subbundle  $W^f$. The image of the tensorial map 
\[
\nabla^f_\cdot R: E\to TF
\]
again lies in $E$ and $R$ is mapped to $\frac{-div^f(R)}{n+1}\cdot R$. Hence $\nabla^fR$ can be interpreted
as a tensorial map on the quotient bundle $Q$. The skew-symmetric part $J$ of this map  
corresponds via $f$ to $J^f$ on $W^f\subset E$. From Lemma \ref{LM8} we know that $J:Q\to Q$
is a complex structure. This complex structure is orthogonal on $Q$ with respect to the 
positive definite scalar product
$f_Q$
which is induced by $f$ on $Q$, i.e., it holds $f_Q(JX,JY)=f_Q(X,Y)$ for all vectors $X,Y$ in $Q$.
Moreover, the $2$-form $f_Q(J\cdot,\cdot)$ equals the $2$-form that is induced by $d\theta$ on the quotient 
$Q$. Thereby, we note that $d\theta(R,\cdot)=0$ on $E$. 

So far we have constructed from $\eR$ with the help of some fixed $f$ in $c$ the data:
\[
\RR R,\quad E,\quad Q:=E/\RR R\quad\mbox{and}\quad J:Q\to Q\ . 
\]
We want to show now that these data do not depend on the choice of $f$, i.e., they are 
uniquely determined by the conformal structure $c$. 
Obviously, the subbundles $\RR R$ and $E:=ker(\theta^f)$ 
do not depend on $f$ in $c$, since $H^1_{\frak{g}}F$ is invariantly identified with $TF$. 
Now, let $\tilde{f}=e^{-2\phi}f$ be an arbitrary rescaled metric in the conformal class $c$. From the
transformation law for Levi-Civita connections in the conformal class we obtain
\[
\nabla_X^{\tilde{f}}R=\nabla_X^fR-d\phi(X) R-d\phi(R)X\qquad \mbox{for\ all}\ X\in E\ .
\]
This relation shows that the skew-symmetric part of $\nabla^{\tilde{f}}R$ induces on
$Q$ the complex structure $J$ again, i.e., $J$ on $Q$ does not depend on
the choice of metric.

As next, we note that, in general, there exists a naturally defined Lie derivative acting on
sections of tractor bundles with respect to conformal Killing vector fields.
This derivative is given by differentiation of a given tractor along the flow of the vector, which 
consists of (local) conformal diffeomorphisms.
For example, let $V$ denote any conformal Killing vector on $F$ and let
$\eQ$ be an adjoint tractor then the Lie derivative of $\eQ$ with respect to $V$ is given by
\[
\eL_V\eQ=\{\eQ,\eS(V)\}\ ,
\]
where $\eS(V)$ is the splitting operator applied to $V$.
In case that $V$ is a Killing vector field for some suitable metric $f$ in $c$ this Lie derivative
coincides with the 
usual $\eL_V$ applied to the components of $\eQ$
with respect to the graded sum
$TF\oplus\frak{co}(TF)\oplus T^*F$.

In particular, for $\eR$ we have
\[
\eL_R\eR=\{\eR,\eR\}=0\ .
\]
In fact, since the conformal vector field $R$ has no zero, we can find locally a metric $f$
in $c$ such that $R$ is Killing, i.e., $R$ is an infinitesimal isometry of $f$ and $\eL_Rf=0$.
The adjoint tractor $\eR$ takes with respect to $f$ the form $(\ R\ ,\ \nabla^f R\ ,\ \eP(f(R,\cdot))\ )$.
The property that $R$ is Killing implies that the commutators $[\eL_R,\nabla^f]$ and $[\eL_R,\eP]$ vanish.
Hence it is 
$[R,R]=0$, $\eL_R(\nabla^fR)=0$ and $\eL_R(\eP(f(R,\cdot)))=0$.
With this remark we can see that $\RR R,E$ and the complex structure 
$J$ on $Q$
are invariant under the (local) flow of the conformal vector field $R$ on $F$. 
In fact, let $\Phi^R_t$ denote this flow with time parameter $t$
and let $f$ be a metric such that $R$ is Killing. Then it holds $[R,R]=0$ and $\eL_R\theta=0$,
which proves that $\Phi^R_{t*}(\RR R)=\RR R$ and $\Phi^R_{t*}(E)=E$
to any time $t$. Moreover,
it is $\eL_R\nabla^fR=0$, which shows that $J$ on $Q$ is invariant under the flow $\Phi^R$ as well.

Furthermore, we notice that the algebraic bracket
\[\begin{array}{cccl}
L_Q :& Q\times Q&\to& K:=TF/E\\[2mm]
&(X,Y)&\mapsto& pr_K[X,Y] \end{array}
\]
on the quotient $Q$ is $\Phi^R$-invariant and non-degenerate. This follows from the 
fact that $E,Q$ are $\Phi^R$-invariant and from the relation 
$\theta^f([X,Y])=-d\theta^f(X,Y)=f(X,\nabla^f_YR)$ for all $X,Y$ in $E$ when $\eL_Rf=0$.
In this situation, it also holds
\[
\theta^f([\nabla^f_XR,\nabla^f_YR])=-d\theta^f(\nabla^f_XR,\nabla^f_YR)=-f(\nabla^f_XR,Y)=f(X,\nabla^f_YR)\ ,
\]
which shows that $L_Q(JX,JY)=L_Q(X,Y)$ for all $X,Y$ in $Q$, i.e., the bracket $L_Q$ is totally real.

Now we consider again the conformal Killing vector $R$ on $F$. 
Since $R$ admits no zero, we can factorise the manifold $F$ locally around every point 
through the integral curves $Int(R)$ of $R$ and obtain a smooth
projection onto a $C^\infty$-manifold of dimension $n$:
\[
\pi_U:\ U\subset F\ \to\ M_U:=U/Int(R)\ .
\]
Thereby, we can assume that the open submanifold $U$ of $F$ is diffeomorphic to
$M_U\times\RR$ and $\pi_U$ is just the natural projection onto the first factor.
Since the subbundle $E$ of $TF$ is invariant under the flow $\Phi^R$, 
it projects by $\pi_U$ to a distribution $H$ in $TM_U$.
The subbundle $\RR R$ of $TF$ vanishes after projection to $M_U$.
Hence the distribution $H$ is of codimension one in $TM_U$ and 
the quotient bundle $Q$ projects naturally onto the distribution $H$.
Since the algebraic bracket $L_Q$ is non-degenerate and $\Phi^R$-invariant on $F$, the distribution $H$ 
is contact in $TM_U$.
Moreover, the $\Phi^R$-invariant complex structure
$J$ on $Q$ projects to a unique complex structure on the contact distribution 
$H$ in $TM_U$, which we again denote by $J$. The algebraic bracket which corresponds to the 
contact distribution $H$ is then real with respect to $J$. 
At this point we have shown that the adjoint tractor $\eR$ on $(F,c)$ 
generates locally a uniquely determined CR-manifold $(M_U,H,J)$ of dimension $n$, which is
partially integrable. 

Eventually, we construct pseudo-Hermitian structures on $(M_U,H,J)$. For this purpose,
let us consider (locally) any metric $f$ in $c$ on $U\subset F$ with $\eL_Rf=0$.
Then, it holds $\eL_R\theta^f=0$ and  $\theta^f$ projects uniquely to a $1$-form
$\theta$ on $M_U$ such that $\theta|_H=0$, i.e., a metric $f$ gives rise to  
a pseudo-Hermitian form $\theta$ on $(M_U,H,J)$ when $\eL_Rf=0$.
From such a $\theta$ we obtain the Levi-form $L_\theta$ on the CR-space
$(M_U,H,J)$. The Levi-form $L_\theta$ is naturally $\pi_U$-related to $f_Q$ on the quotient bundle $Q$,
which shows that $(M_U,H,J)$ and $\theta$ are strictly pseudoconvex. 
\begin{PR} \label{PR4} Let $(F,c)$ be a space with conformal structure $c$ admitting
$\eR\in\Gamma(\mathcal{A}(F))$ which acts as complex structure
on $\eT(F)$ such that the 'first slot' $R$ is conformal Killing. Then
\begin{enumerate}
\item
the local factorisation of $F$ through the integral curves $Int(R)$ to $R$ 
generates a smooth space admitting a uniquely determined strictly pseudoconvex and partially integrable CR-structure. 
\item 
Any (local) metric $f$ in $c$ with $\eL_Rf=0$ generates a pseudo-Hermitian
structure on that CR-space.
\end{enumerate}
\end{PR}

So far we have not used in our discussion all the information that is encoded in 
the existence of the complex structure $\eR$. And, in fact, there is still an open issue.
It is
the question how the conformal structure $c$ on $U\subset F$ is related to the induced 
CR-space $(M_U,H,J)$. To give an answer to this question we consider the $1$-form 
$A^f:=\eP(f(R,\cdot))$ with respect to any metric $f$ such that $\eL_Rf=0$.
Then it holds $A^f(R)=-1$ and $\eL_RA^f=0$. These properties show that we can understand
$A^f$ as a 'local connection form' on the fibration $\pi_U:U\subset F\to M_U$. 
From our construction so far, it is clear that $f$ on $U$ takes the form
\[
\pi_U^*L_\theta-2\cdot\theta^f\circ A^f\ ,
\]
where $L_\theta$ is the Levi-form on $(M_U,H,J)$ belonging to $\theta$, which in turn
was induced by $f$ on $U\subset F$.

If we choose $U$ suitably small we can identify it with an open subset of the canonical $S^1$-bundle 
over $(M_U,H,J)$. (Such an identification can be understood as the choice of a
gauge in the canonical $S^1$-bundle.) 
With such a gauge we can compare $f$ on $U$ with the Fefferman construction 
to $\theta$ on the canonical $S^1$-bundle of $(M_U,H,J)$. So let $A_\theta$ be
the Weyl connection to $(M_U,H,J)$ that we introduced in section \ref{ab5}. We denote by
$\ell:=A^f-A_\theta$ the difference of the two given connections (via the chosen gauge).
Remember that $A_\theta$ is determined by $f$. We can see that $f=L_\theta-2\cdot\theta^f\circ A^f$ is just 
the $\ell$-Fefferman metric on 
$U$ over $(M_U,H,J)$. However, the difference $\ell$ between the two connections can not be arbitrary. In fact,
the Killing vector $R$ is vertical with respect to $\pi_U$, i.e., $R$ is vertical 
in the $\ell$-Fefferman construction
when we consider $U$ as a subset of the canonical $S^1$-bundle over $M_U$. By our assumption, $\eS(R)=\eR$ 
is a complex structure. Yet we know from Proposition \ref{PR3} that a vertical Killing 
vector in the $\ell$-Fefferman construction gives rise to a complex structure via the splitting operator $\eS$ 
if and only if the contraction $tr_\theta d\ell(\cdot,J\cdot)=0$ on $M_U$ vanishes, 
i.e., (the $\pi_U$-projection of) 
$\ell$ is an 
element of $H^1_{tr}(M_U,H,J)$. This completes our discussion of the general reconstruction.

\begin{PR} \label{PR5} Let $(F,c)$ be a space with conformal structure $c$ admitting
$\eR\in\Gamma(\mathcal{A}(F))$ which acts as complex structure
on $\eT(F)$ such that the 'first slot' $R$ is conformal Killing
and let $(M_U,H,J)$ be the (locally) induced CR-space (cf. Proposition \ref{PR4}).
Then $c$ on $F$ is (locally) conformally
equivalent to some $\ell$-Fefferman metric which is constructed on $(M_U,H,J)$
with suitable $[\ell]\in H^1_{tr}(M_U,H,J)$.
\end{PR}
We remark that the Nijenhuis tensor $\eN_J$ of the CR-space $(M_U,H,J)$ and the $1$-form class $[\ell]$
can be recovered from the curvature expression $\Omega^{nor}(R,\cdot)$.

It is well known that
the classical Fefferman construction over integrable CR-spaces produces conformal classes $c$ of metrics,
which admit a parallel complex structure on the standard tractor bundle $\eT(F_c)$, i.e., 
$\nabla^{nor}J_{CR}=0$  for $J_{CR}=\eS(R)$ and $\Omega^{nor}(R,\cdot)=0$. The latter condition is equivalent
to $R\ecke \eW=0$ and $R\ecke \eC=0$, where $\eW$ denotes the Weyl tensor of the conformal structure $c$
and $\eC$ denotes the Cotton-York tensor. The Weyl tensor $\eW$ is determined by the $\frak{g}_0$-part of the 
curvature $\Omega^{nor}$ and $\eC$ comes from the $\frak{g}_{1}$-part. 
On the other side, the existence of a lightlike Killing vector $R$ for some Lorentzian metric
$f$ satisfying
\[
R\ecke \eW=0,\qquad R\ecke \eC=0\qquad \mbox{and}\qquad Ric(R,R)>0
\]
implies that $f$ is locally isometric to the Fefferman metric of an integrable CR-space (cf. \cite{Gra87}).
In the situation when there exists a
complex structure $J_{CR}$ with $\nabla^{nor}J_{CR}=0$, these conditions are automatically satisfied
with respect to a metric $f$ such that the 'first slot' of $J_{CR}$ is (locally) a Killing vector. 
The reason is that $\eP(f(R,\cdot))(R)=-1$ implies $Ric(R,R)>0$ (cf. section \ref{ab9}). 
Notice once again that such a metric 
$f$ always 
exists
locally.

In particular, this shows that the existence of a parallel complex structure $J_{CR}$ on 
the standard tractor bundle $\eT(F)$ gives rise at least locally to a parallel section in
the canonical complex line tractor bundle, which is defined as  \[\eF:=\Lambda^{m+2*}(\eT_{10}),\] where
$\eT_{10}$ denotes the $i$-eigenspace bundle to the extended endomorphism $J_{CR}$ on the complexified
standard tractor bundle $\eT^\CC(F)=\eT(F)\otimes\CC$. In fact, 
the conformal curvature $\Omega^{nor}$ acts on $\eT(F)$ by
\[
\Omega^{nor}(X,Y)(t)=\kappa(X,Y)\bullet t
\]
for all $t\in\eT(F)$ and vectors $X,Y\in\ TF\cong H^1_\frak{g}F$. 
Now let $f$ be any metric in $c$, let $U$ denote the $f$-orthogonal complement 
to the span of $R=\pi_H(J_{CR})$ and the dual vector field of $\eP(f(R,\cdot))$
and let $J$ by the antisymmetric part of the endomorphism $\nabla^fR:TF\to TF$.
Furthermore, let $\{e_i: i=1,\ldots,2m\}$ denote an orthonormal basis of $(U,f|_{U})$
such that 
$J(e_{2\alpha-1})=e_{2\alpha}$ and $J(e_{2\alpha})=-e_{2\alpha-1}$ for all 
$\alpha=1,\ldots,m$. Then it is a straightforward calculation to see that the conformal curvature 
$\Omega^{nor}$ acts on the line tractor bundle $\eF$ by multiplication with the complex number $\rho$,
which is given by
\[
i\cdot\sum_{\alpha=1}^{m}e^*_{2\alpha}(\kappa_0(X,Y)(e_{2\alpha-1}))=
i\cdot\sum_{\alpha=1}^{m}\eW(X,Y,e_{2\alpha-1},Je_{2\alpha-1})\ .\]
The property of $J_{CR}$ to be parallel for $\nabla^{nor}$ implies
directly that \[\eW(X,Y,JZ,V)=-\eW(X,Y,Z,JV)\qquad \mbox{for\ all}\ \ X,Y,Z,V\in TF .\]
Using the Bianchi identity and the fact that $\partial^*\circ\kappa=0$  
(i.e. the fact that the Weyl tensor $\eW$ has no trace) 
we find that the complex number $\rho$ is zero, i.e.,
the line tractor bundle $\eF$ has no conformal curvature. Hence $\eF$ admits locally
parallel sections. 
\begin{THEO} \label{THex}(cf. \cite{Spa85}, \cite{Gra87}) Let $(F,c)$ be a manifold with a conformal
structure $c$ of Lorentzian signature admitting a conformal Killing vector $R$
whose adjoint tractor $\eS(R)$ is a $\nabla^{nor}$-parallel complex structure. Then $c$ is locally conformally
equivalent to a Fefferman metric $f_\theta$ constructed over an integrable CR-space.
\end{THEO}
In particular, Theorem \ref{THex} says that \[\Omega^{nor}(R,\cdot)=0\qquad \mbox{if\ and\ only\ if}\qquad 
\eN_J=0\ \ \mbox{and}\ \ [\ell]=0\ ,\]
where $R$ is a conformal Killing vector and $\eS(R)$ is a complex structure. 
In terms of conformal holonomy the statement of Theorem \ref{THex} says that if the holonomy algebra of
the normal conformal connection $\omega^{nor}$ is reduced to $\frak{u}(1,m)$ then it is
already reduced to $\frak{su}(1,m)$. The latter is the Lie algebra of the structure group 
$\SU(1,m)$ of CR-geometry. In particular, if a complex structure $J_{CR}$ with $\nabla^{nor}J_{CR}=0$ 
exists, then there exists locally also a solution of the twistor equation for spinors (cf. \cite{Bau99}).

\section{Main Theorem}
\label{ab12}

Proposition \ref{PR3} and \ref{PR5} establish  our main result about 
complex structures on the conformal standard tractor bundle. We combine them
in order to formulate the following Theorem. We remember that $c_{[\ell]}$ denotes
the conformal class which is determined by a $\ell$-Fefferman metric $f_{\theta,\ell}$,
where $[\ell]$ denotes the class of a $1$-form $\ell$ modulo closed forms.
  
\begin{THEO} \label{TH3} Any smooth manifold $(F^{n+1},c)$ with Lorentzian conformal structure $c$
admitting a conformal Killing vector field $R$ whose corresponding adjoint tractor
$\eS(R)$ is a complex structure on the standard tractor bundle $\eT(F)$ 
is locally conformally equivalent to a $[\ell]$-Fefferman conformal class constructed
on some  
strictly pseudoconvex, partially integrable CR-structure $(M,T_{10})$
with some $[\ell]\in H^1_{tr}(M,T_{10})$, i.e.,  the map
\[
\Psi:\ (\ T_{10}\ ,\ [\ell]\ )\qquad\mapsto\qquad (\ c_{[\ell]}\ , \ \eS(R)\ )
\]
is surjective for $[\ell]\in H^1_{tr}(M,T_{10})$ onto the local conformal classes with complex structures
$\eS(R)$.
\end{THEO}   
We remark 
that it is not clear at this point whether the map $\Psi$ is a bijection. In principle,
there might be the possibility that a conformal class is $[l]$-Fefferman with respect to two different
constructions along different fibres (resp. lightlike Killing vector fields), which are not
equivalent. 
We remark also that
non-trivial examples for elements in $H^1_{tr}(M,T_{10})$ do exist for CR-structures, in
general. However, we do not aim to investigate here the space $H^1_{tr}(M,T_{10})$ and the 
properties of their elements.  
We expect that the $\ell$-Fefferman construction is really an extension of the Fefferman construction
and produces new sorts of Lorentzian metrics and conformal classes.

There are several more aspects of the relation of CR-geometry and conformal geometry, 
which are not touched here. We want to mention some of them. 
This happens also with
the intention to set our discussion into a broader context. In fact, an issue that is behind 
our investigations here, is the relation of the canonical normal Cartan connections of CR-geometry
and conformal geometry in the Fefferman construction. 
Since the Fefferman construction can also be generalised and applied to 
various other parabolic geometries, this topic seems to be very basic and of general importance.

For integrable CR-structures it is well known that the equivariant extension of the normal 
CR-Cartan connection in the Fefferman construction gives rise to the normal Cartan connection
of conformal geometry. In fact, this is true if and only if the CR-geometry is integrable.
In case that the CR-geometry is only partially integrable, there still exists a unique 
normal Cartan connection (cf. \cite{CS00}). However, the equivariant extension in the Fefferman construction 
does not give rise to the normal connection of conformal geometry. 

What we have shown in this paper is that the two natural connections in the Fefferman construction
(without an $\ell$-term)
of partially integrable CR-structures 
coincide at least when they are applied to the natural
complex structure $J_{CR}$ that appears in the construction. Of course, this does not mean 
that they are the same. In fact, it shows only that they are equal up to a $1$-form with values
in the Lie algebra $\frak{u}(1,m)$.
It remains to determine this 
$1$-form 
in the difference of the two natural Cartan connections. To achieve this, it should be in principle enough
to calculate all components of the Ricci curvature of the Fefferman metric with respect to some pseudo-Hermitian
structure.    

Moreover, if we know the Ricci curvature of the Fefferman metric we can proceed and calculate the conformal
curvature in dependence of the Tanaka-Webster curvature. So we would know then the right hand side
$-\Omega^{nor}(\pi_H(\eR),\cdot)$ of our tractor equation in terms of Webster curvature. This 
curvature expression would characterise 
essential properties
of the fundamental Killing vector in the Fefferman construction. However, one reason that 
we did not continue to go this way was that calculations are already rather extensive. 
And, since the Fefferman construction is invariant in its nature, the idea is that there should be 
an invariant approach and calculation for the comparison of the two normal connection. 
This means an approach that does not use a Weyl connection or any other form of preferred 
linear connection. We do not know yet how this invariant approach works. On the other side, 
the occurrence of the $\ell$-Fefferman construction seems to show that the existence of complex
structures on standard tractors 
should not only be seen as a CR-invariant problem. Maybe this is an indication that an invariant approach
to the comparison of the normal Cartan connections is not too obvious in the first place.

\section*{Acknowledgements}
This work was done during a stay at the Erwin Schr{\"o}dinger International (ESI)
Institute for Mathematical Physics in Vienna/Austria.
I would like to thank Andreas Cap for supporting me with the research to this work
through many helpful discussions.
My stay at the ESI was funded by a Junior Research Fellowship grant of
the Austrian Ministry of Education, Science and Culture (BMBWK).


\end{sloppypar}
\end{document}